\documentclass{amsart} 
\usepackage{amssymb,latexsym}

\newcommand{\Pb}{\operatorname{\bf P}} 
\newcommand{\Qb}{\operatorname{\bf Q}} 
\newcommand{\Rb}{\operatorname{\bf R}} 
 
\newcommand{\ztj}{\mathbb{Z}_{t_j}} 
\newcommand{\mult}{{\rm mult}} 
 
\newcommand{\Q}{\mathbb{Q}} 
\newcommand{\QQ}{\mathcal{Q}} 
\newcommand{\C}{\mathcal{C}} 
\newcommand{\F}{\mathcal{F}} 
\newcommand{\M}{\mathcal{M}} 
 
\newcommand{\K}{\mathcal{K}}

\theoremstyle{plain} 
\newtheorem{theorem}{Theorem}[section] 
\newtheorem{corollary}[theorem]{Corollary} 
\newtheorem{lemma}[theorem]{Lemma} 
\newtheorem{proposition}[theorem]{Proposition} 
\newtheorem{definition}[theorem]{Definition} 
 
\theoremstyle{remark} 
\newtheorem{remark}{Remark}[section] 
\newtheorem{example}[remark]{Example} 
\newtheorem{conjecture}[remark] {Conjecture} 
 
\numberwithin{equation}{section}

\begin{document} 
\title[Decomposable Compositions] 
{Decomposable Compositions, Symmetric Quasisymmetric 
Functions and Equality of Ribbon Schur Functions} 
\author{Louis J.\ Billera} 
\address{Department of Mathematics, Cornell University, Ithaca, NY 
14853-4201 USA} 
\email{billera@math.cornell.edu} 
\thanks{The first author was supported in part by NSF grant 
DMS-0100323.  The second and third authors were supported in part by the 
National Sciences and Engineering Research Council of Canada.} 
\author{Hugh Thomas} 
\address{Fields Institute, 222 College St., Toronto, ON M5T 3J1 Canada} 
\curraddr{Department of Mathematics and Statistics, University of New 
Brunswick, Fredericton, NB E3B 5A3 Canada} 
\email{hugh@math.unb.ca} 
\author{Stephanie van Willigenburg} 
\address{Department of Mathematics, University 
of British Columbia, Vancouver, BC  V6T~1Z2 Canada } 
\email{steph@math.ubc.ca} 
 
 
\begin{abstract} 
We define an equivalence relation on integer compositions and 
show that two ribbon Schur functions are identical if and only 
if their defining compositions are equivalent in this sense. This 
equivalence is completely determined by means of a 
factorization for compositions: equivalent compositions have 
factorizations that differ only by reversing some of the terms. 
As an application, we can derive identities on certain 
Littlewood-Richardson coefficients. 
 
Finally, we consider the cone of symmetric functions having a 
nonnnegative representation in terms of the fundamental 
quasisymmetric basis.  We show the Schur functions are among the 
extremes of this cone and conjecture its facets are 
in bijection with the equivalence classes of compositions. 
\end{abstract} 
 
\maketitle

%
 
\section{Introduction \label{intro}} 
 
An important basis for the space of symmetric functions of degree $n$ is the 
set of  classical Schur functions  $s_\lambda$, where $\lambda$ runs over all 
\emph{partitions} of $n$. 
The skew Schur functions $s_{\lambda/\mu}$ can be 
expressed in terms of these by means of the Littlewood-Richardson 
coefficients $c^\lambda_{\mu\nu}$ by 
\begin{equation} 
s_{\lambda/\mu}=\sum_{\nu} c^\lambda_{\mu\nu} s_\nu. 
\label{LRcoeff} 
\end{equation} 
These coefficients also describe the structure constants in the 
algebra of symmetric functions. In particular they describe the 
multiplication rule for Schur functions, 
\begin{equation} 
s_\mu \ s_\nu = \sum_\lambda c^\lambda_{\mu\nu} s_\lambda. 
\end{equation} 
{}From the perspective of the representation theory of the 
symmetric group, the coefficient $c^\lambda_{\mu\nu}$ gives the 
multiplicity of the irreducible 
representation corresponding to the partition $\lambda$ in the 
induced tensor product of those corresponding to $\mu$ and $\nu$.  In 
algebraic geometry the $ c^\lambda_{\mu\nu}$ arise as intersection 
numbers in the Schubert calculus on a Grassmanian. As a result of 
these and other instances in which they arise, the determination 
of these coefficients is a central problem. 
 
We consider here the question of when two ribbon Schur 
functions might be equal.  The more general question of equalities among 
all skew Schur functions has since been broached in \cite{SvW}, where the 
question of when a skew Schur function can equal a Schur function is answered  
in the case of power series and 
their 
associated polynomials. 
The general question  
of equalities among skew Schur functions remains open.  Any such equalities 
imply equalities between certain pairs of Littlewood-Richardson coefficients.

In the case 
of \emph{ribbon} Schur functions, that is skew Schur functions 
indexed by a shape known as a ribbon (or rim hook, or border 
strip; these are diagrams corresponding to connected skew shapes containing 
no $2\times2$ rectangle), we give necessary and sufficient conditions for equality. 
Ribbons are in natural correspondence with compositions, and 
equality arises from an equivalence relation on compositions, whose 
equivalence classes all have size equal to a power of two. This 
power corresponds to the number of nonsymmetric compositions in a 
certain factorization of any of the underlying compositions in a 
class, and equivalence comes by means of reversal of terms. 
 
A motivation for studying ribbon Schur functions is that they 
arise in various contexts.   They were studied already by MacMahon 
\cite[\S 168-9]{MM}, who showed their coefficients in terms of the 
monomial symmetric functions to count descents in permutations 
with repeated elements. 
The scalar product of any two gives the 
number of permutations such that it and its inverse have the 
associated pair of descent sets \cite[Corollary 7.23.8]{ECII}. 
They are also useful in computing the number of permutations with 
a given cycle structure and descent set \cite{GR}.  Lascoux and 
Pragacz \cite{LP} give a determinant formula for computing Schur 
functions from associated ribbon Schur functions.  
In addition,  
they arise  as $sl_{n}$-characters of the irreducible components 
of the Yangian representations in level 1 modules of $\widehat{sl}_{n}$ 
\cite{KKN}. 
 
In the theory of noncommutative symmetric functions of Gel'fand  
{\it et al.\ }\cite{GKL},  
the noncommutative analogues of the ribbon Schur functions 
form a homogeneous linear basis.  It is therefore of some interest to know what 
relations are introduced when passing to the commutative case.  In 
particular, which pairs become identical? In  \cite[Exercise 7.56 b)]{ECII} the ribbon Schur function indexed by a composition and its reversal are seen to be identical.  However, as we shall see, this is not the whole story. 
 
This paper is organized as follows. In Section 
\ref{schur-section}, we introduce an equivalence relation on 
compositions and derive some if its properties.  The relation is 
defined in terms of coefficients of symmetric functions when 
expressed in terms of the fundamental basis of the algebra of 
quasisymmetric functions.  We show that this relation can be 
viewed combinatorially in terms of  coarsenings of the respective 
compositions. 
Theorem \ref{ribbon-equality} then 
shows compositions to be equivalent if and only if their 
corresponding ribbon Schur functions are identical. 
 
Section \ref{composition-sec} introduces a binary operation on 
compositions.  In the case of compositions denoting the descent 
sets of a pair of permutations, the operation results in the 
composition giving the descent set  of their tensor product.  In 
Sections \ref{decomp-proof} and \ref{technical} we prove our main 
result, Theorem \ref{decomp-thm}, which states that equivalence of 
two compositions is precisely given by reversal of some or all of 
the terms in some factorization. Thus the congruence classes all 
have size given by a power of two; this power is the number of 
nonsymmetric terms in the finest factorization of any composition 
in this class. 
 
Finally, in Section \ref{F-pos}, we consider the cone of $F$-positive symmetric 
functions, showing the Schur functions to be among its extremes and 
conjecturing 
its facets to be in one-to-one correspondence with equivalence classes of 
compositions. 
 
The remainder of this section contains the basic definitions we will be using.  Where 
possible, we are using the notation of \cite{MacD} or \cite{ECII}. 
 
\subsection{Partitions and compositions.} A composition $\beta$ of $n$, denoted $\beta\vDash n$, is a 
 list of positive integers $\beta _1\beta _2\ldots\beta _k$ such that $\beta _1+\beta _2+\ldots +\beta_k=n$. 
 We refer to each of the $\beta _i$ as components, and say that 
$\beta$ has \emph{length} 
 $l(\beta)=k$ and \emph{size} $|\beta |=n$. If the components of $\beta$ are weakly decreasing 
 we call $\beta$ a \emph{partition}, denoted $\beta \vdash n$ and refer to each of the $\beta _i$ as parts. 
 For any composition $\beta$ there will be two other closely 
 related compositions that will be of interest to us. The first is the reversal of $\beta$, 
 $\beta ^\ast = \beta_k\ldots\beta _2\beta_1$, and the second is the partition determined by $\beta$, 
 $\lambda(\beta)$, which is obtained by reordering the components of $\beta$ in weakly decreasing 
 order, e.g. $\lambda(3243)=4332$. Moreover we say two compositions $\beta, \gamma$ determine the 
  same partition if $\lambda(\beta)=\lambda(\gamma)$. 
 
Any composition $\beta \vDash n$ also naturally corresponds to a subset $S(\beta)\subseteq [n-1] 
=\{ 1,2,\ldots , n-1\}$ where $$S(\beta)=\{\beta _1,\beta _1+\beta _2,\beta _1+\beta _2+\beta _3, 
 \ldots  ,\beta_1+\beta _2 +\ldots +\beta _{k-1}\}.$$ Similarly any subset $S=\{i_1,i_2\ldots , 
  i_{k-1}\}\subseteq [n-1]$ corresponds to a composition $\beta (S)\vDash n$ where 
  $$\beta(S)=i_1 (i_2-i_1) (i_3-i_2)\ldots (n-i_{k-1}).$$

Finally, recall two partial orders that exist on compositions. We say that for compositions 
$\beta,\gamma \vDash n$, we write $\beta \prec \gamma$ when $\beta$ is 
 \emph{lexicographically less} than $\gamma$, that is, 
$\beta = \beta_1\beta_2\cdots \neq 
 \gamma_1\gamma_2\cdots = \gamma$, and the first $i$ for which $\beta_i \neq \gamma_i$ 
  satisfies 
$\beta_i < \gamma_i$.  In particular, $11\cdots 1 \preceq \beta \preceq n$ for any 
  $\beta \vDash n$. Secondly, given any two compositions $\beta$ and $\gamma$ we say 
$\beta$ is a \emph{coarsening} of $\gamma$, denoted $\beta\ge\gamma$, if we can obtain 
$\beta$ by adding together adjacent components of $\gamma$, e.g., $3242\ge 3212111$. 
 Equivalently, we can say $\gamma$ is a \emph{refinement} of $\beta$. 
 
\subsection{Quasisymmetric and symmetric functions.} We denote by $\QQ$ the algebra of 
quasisymmetric functions over $\Q$, that is all bounded degree 
formal power series $F$ in variables $x_1,x_2,\ldots $ such that 
for all $k$ and $i_1<i_2<\ldots <i_k$ the coefficient of 
$x_{i_1}^{\beta_1}x_{i_2}^{\beta _2}\ldots x_{i_k}^{\beta _k}$ is 
equal to that of $x_1^{\beta_1}x_2^{\beta _2}\ldots x_k^{\beta 
_k}$. There are two natural bases for $\QQ$ both indexed  by 
compositions $\beta =\beta _1\beta _2\ldots\beta _k$, $\beta_i>0$: 
 the monomial basis spanned by $M_0=1$ and all power series $M_\beta$ where 
$$M_\beta =\sum _{ i_1<i_2<\ldots<i_k}x_{i_1}^{\beta _1}x_{i_2}^{\beta _2}\ldots x_{i_k}^{\beta _k}$$ 
and the fundamental basis 
spanned by $F_0=1$ and all power series $F_\beta$ where 
$$F_\beta = \sum _{\gamma\leq \beta }M_\gamma .$$ 
Note that $\QQ$ is a graded algebra, with $\QQ_n = \text{span}_\Q \{ M_\beta ~|~ \beta \models n \}$. 
 
We define the algebra of symmetric functions $\Lambda$ to be the subalgebra of $\QQ$ 
spanned by the \emph{monomial symmetric functions} 
\begin{eqnarray} 
m_\lambda = \sum_{\beta:\lambda(\beta)=\lambda} M_\beta,\quad 
\lambda\vdash n, ~n>0 \label{monomial} 
\end{eqnarray} 
and $m_0=1$. Again, $\Lambda$ is graded, with $\Lambda_n = \Lambda 
\cap \QQ_n$. 
 
From quasisymmetric functions we can define Schur functions, which also form a basis 
for the symmetric functions, but first we need to recall some facts about  tableaux. 
 
For any partition $\lambda = \lambda _1\ldots\lambda _k\vdash n$ the related \emph{Ferrers diagram} 
 (by abuse of notation also referred to as $\lambda$) is an array of left justified boxes with $\lambda _1$ 
  boxes in the first row, $\lambda _2$ boxes in the second row, and so on.  For example, 
the Ferrers diagram 4332 is 
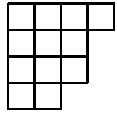
\begin{center} 
\begin{picture}(100,50)(-30,-5) 
\put(0,0){\line(1,0){20}} 
\put(0,10){\line(1,0){30}} 
\put(0,20){\line(1,0){30}} 
\put(0,30){\line(1,0){40}} 
\put(0,40){\line(1,0){40}} 
 
\put(0,40){\line(0,-1){40}} 
\put(10,40){\line(0,-1){40}} 
\put(20,40){\line(0,-1){40}} 
\put(30,40){\line(0,-1){30}} 
\put(40,40){\line(0,-1){10}} 
\end{picture} 
\end{center} 
 
A \emph{(Young) tableau} of shape $\lambda$ and size $n$ is a filling of the boxes of 
$\lambda$ with positive integers. If the rows weakly increase and the columns strictly increase we 
say it is a \emph{semi-standard} tableau, and if in addition, the filling of the boxes involves the integers 
$1, 2, \ldots ,n$ appearing once and only once we say it is a \emph{standard} tableau. Note that in this 
instance  both the rows and columns strictly increase.  
  
 More generally, we can define skew diagrams and skew tableaux.  Let $\lambda , \mu$ be 
partitions such that if there is a box in the $(i,j)$-th position in the Ferrers diagram $\mu$ then there is a 
box in the $(i,j)$-th position in the Ferrers diagram $\lambda$. The skew diagram $\lambda /\mu$ is the 
array of boxes $\{c\ |\ c\in\lambda , \ c\not\in\mu\}$.  
For example, the skew diagram 4332/221 is 
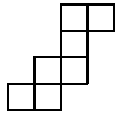
\begin{center}\begin{picture}(100,50)(-30,-5) 
 
\put(0,0){\line(1,0){20}} 
\put(0,10){\line(1,0){30}} 
\put(10,20){\line(1,0){20}} 
\put(20,30){\line(1,0){20}} 
\put(20,40){\line(1,0){20}} 
 
\put(0,10){\line(0,-1){10}} 
\put(10,20){\line(0,-1){20}} 
\put(20,40){\line(0,-1){40}} 
\put(30,40){\line(0,-1){30}} 
\put(40,40){\line(0,-1){10}}\end{picture}\end{center} 
We can then define skew tableaux, semi-standard 
skew tableaux, and standard skew tableaux analogously.   
  
 Given a standard  tableau or skew tableau $T$, we say it has a 
 \emph{descent} in position $i$ if $i+1$ appears in a lower row than $i$. 
 Denote the set of all descents of $T$ by $D(T)$.   We take \cite[Theorem 7.19.7]{ECII} as 
 our definition of skew Schur functions.

 
 
\begin{definition}\label{schurasf} 
Let $\lambda, \mu \vdash n$ such that $\lambda/\mu$ is 
a skew diagram.  Then the skew Schur function $s_{\lambda/\mu}$ is defined by 
\begin{eqnarray} 
\label{schur} 
s_{\lambda/\mu} = \sum_T F_{\beta(D(T))} 
\end{eqnarray} where the sum is over all standard tableaux $T$ of shape $\lambda/\mu$. 
\end{definition} 
 
The \emph{Schur functions} $s_{\lambda}$ are those skew Schur functions with $\mu = 0$. 
For example, $s_{22}= F_{22}+F_{121}$.

A skew diagram is said to be \emph{connected} if, regarded as a 
union of squares, it has a connected interior. If the skew diagram 
$\lambda/\mu$ is connected and contains no $2\times 2$ array of 
boxes we call it a \emph{ribbon}. Observe ribbons of size $n$ are 
in one-to-one correspondence with compositions $\beta$ of size $n$ 
by setting $\beta _i$ equal to the number of boxes in the $i$-th 
row from the bottom.  For example, the skew diagram 4332/221 
%
%
is a ribbon, corresponding to the composition 2212. 
 
Henceforth, we will denote ribbons by compositions, and denote the 
skew Schur functions $s_{\lambda/\mu}$, for a ribbon  $\lambda/\mu$, 
as the \emph{ribbon Schur function} $r_{\beta}$, where $\beta$ is the 
corresponding composition.  Thus $r_{2212} := s_{4332/221}$. 
 
 
Further details on symmetric functions can be found in \cite{ECII}. 
 
\subsection{Acknowledgements.}  The authors are grateful to 
Marcelo Aguiar, Nantel Bergeron, Curtis Greene, Alexander Postnikov,
and Richard Stanley for useful comments over the course of
the work,  and to GAP \cite{Gap} for helpful calculations.

\section{Equality of Ribbon Schur Functions}\label{schur-section} 
 
Although, in general, it is difficult to determine when two skew 
Schur functions are equal, it transpires that when computing 
ribbon Schur functions equality is determined via a 
straightforward equivalence on compositions. 
 
\subsection{Relations on ribbon Schur functions} 
 
There is a useful representation of ribbon Schur functions in terms of the basis of complete homogeneous 
symmetric functions $h_{\lambda}=h_{\lambda_{1}} h_{\lambda_{2}} \cdots h_{\lambda_{k}}$, 
known already to MacMahon \cite[\S 168]{MM}.  It can be derived from the Jacobi-Trudi identity 
\cite[Theorem 7.16.1]{ECII} 
\begin{equation}\label{jacobitrudi} 
s_{\lambda/\mu}= \det (h_{\lambda_{i}-\mu_{j}-i+j}), 
\end{equation} 
where $h_{0}=1$ and $h_{k}=0$ if $k<0$. 
 
\begin{proposition} 
\label{ribbonform} 
For any $\alpha \vDash n$, 
\begin{equation*} 
r_{\alpha}=  (-1)^{l(\alpha)} \sum_{\beta \ge \alpha}(-1)^{l(\beta)} h_{\lambda(\beta)}. 
\end{equation*} 
\end{proposition} 
\begin{proof} 
Applying (\ref{jacobitrudi}) to the ribbon shape $\lambda/\mu$ corresponding to $\alpha=\alpha_{1}\cdots \alpha_{k}$, 
we get 
\begin{equation*} 
r_{\alpha} = \det \left[ 
\begin{matrix} 
h_{\alpha_{k}}&h_{\alpha_{k-1}+\alpha_{k}}&h_{\alpha_{k-2 }+\alpha_{k-1}+\alpha_{k}}&\cdots& 
h_{\alpha_{1}+\cdots+\alpha_{k}} \cr 
1& h_{\alpha_{k-1}} & h_{\alpha_{k-2}+\alpha_{k-1}} & \cdots & h_{\alpha_{1}+\cdots+\alpha_{k-1}} \cr 
&1&h_{\alpha_{k-2}}&\cdots& h_{\alpha_{1}+\cdots+\alpha_{k-2}} \cr 
&&\ddots  &&\vdots \cr 
&&&1&h_{\alpha_{1}} 
\end{matrix} \right]. 
\end{equation*} 
Expanding down the first column gives 
$$r_{\alpha}= r_{\alpha_{1}\cdots \alpha_{k-1}}h_{\alpha_k} 
- r_{\alpha_{1}\cdots \alpha_{k-1}+\alpha_k},$$ 
which along with induction on $l(\alpha)$ gives the desired result. 
\end{proof} 
 
By inverting the relation in Proposition \ref{ribbonform}, we see immediately that the 
ribbon Schur functions $r_{\alpha}$ generate $\Lambda$. 
 
It is straightforward to establish all the algebraic relations that hold among 
ribbon Schur functions.  Using Proposition \ref{ribbonform}, one can show 
that ribbon Schur functions satisfy the multiplicative relations 
\begin{equation}\label{ribbonmult} 
r_{\alpha}\thinspace r_{\beta} = r_{\alpha \cdot \beta} + r_{\alpha \odot \beta}, 
\end{equation} 
where for $\alpha =\alpha_{1}\cdots \alpha_{k}$ and $\beta= \beta_{1}\cdots \beta_{l}$, 
$ \alpha \cdot \beta = \alpha _1\ldots \alpha _k\beta _1\ldots 
\beta _l$ 
is the usual operation of \emph{concatenation} 
and 
$ \alpha\odot\beta = \alpha _1\ldots \alpha _{k-1}(\alpha _k + 
\beta _1)\beta _2 \ldots \beta _l$, 
is  \emph{near concatenation}, 
which differs from concatenation in that the last component of 
$\alpha$ is added to the first component of $\beta$.   
The relation (\ref{ribbonmult})  has been known since MacMahon \cite[\S 169]{MM}. 
Proofs of (\ref{ribbonmult}) in the noncommutative setting can be found in 
\cite[Prop. 3.13]{GKL} and \cite[Prop. 4.1]{BLiu}. 
We show next that these relations generate all the relations 
among ribbon Schur functions. 
 
\begin{proposition} 
\label{ribbonrelations} 
Let {$z_{\alpha}, \alpha \vDash n, n\ge 1$} be commuting indeterminates.  Then 
as algebras, ${\Lambda}$ is isomorphic to the quotient 
$${\Q[z_{\alpha}]/ 
 \langle z_{\alpha}\ z_{\beta}-z_{\alpha\cdot \beta} - z_{\alpha \odot \beta}\rangle}.$$ 
\end{proposition} 
\begin{proof} 
Consider the map $\varphi:\Q[z_{\alpha}] \rightarrow \Lambda$ defined by 
 $z_{\alpha} \mapsto r_{\alpha}$.  This map is surjective since the $r_{\alpha}$ 
generate $\Lambda$.   Grading $\Q[z_{\alpha}]$ by setting the degree of $z_{\alpha}$ 
to be $n=|\alpha|$ makes $\varphi$ homogeneous. 
To see that $\varphi$ induces an isomorphism with the quotient, note that 
$\Q[z_{\alpha}]/ 
 \langle z_{\alpha}\ z_{\beta}-z_{\alpha\cdot \beta} - z_{\alpha \odot \beta}\rangle$ 
 maps onto $\Q[z_{\alpha}]/\ker \varphi \simeq\Lambda $, 
 since $ \langle z_{\alpha}\ z_{\beta}-z_{\alpha\cdot \beta} - z_{\alpha \odot \beta}\rangle 
 \subset \ker\varphi$. 
  
 It will suffice to show that  the degree $n$ component of 
$\Q[z_{\alpha}]/ 
 \langle z_{\alpha}\ z_{\beta}-z_{\alpha\cdot \beta} - z_{\alpha \odot \beta}\rangle$  is generated by 
the images of the $z_{\lambda}, \lambda \vdash n$, and so has dimension at 
most the number of partitions of $n$.  
We have the relations 
\begin{equation}\label{ribbonlinear} 
z_{\alpha \cdot \beta}+z_{\alpha\odot \beta} = 
z_{\beta \cdot \alpha} + z_{\beta \odot \alpha} 
\end{equation} 
by commutativity of the $z_{\alpha}$. 
Let $\gamma =g_1\dots g_k\vDash n$.  We will show by induction on
$k$ that $z_\gamma=z_{\lambda(\gamma)} +$ 
a sum of $z_\delta$'s with $\delta$ having no more than $k-1$ parts.

Let $g_i$ be 
a maximal component of $\gamma$.  Then by (\ref{ribbonmult}) and
(\ref{ribbonlinear}) we have
\begin{eqnarray*}
z_\gamma&=&
z_{g_i\dots g_k g_1\dots g_{i-1}} + 
\textrm{a sum of $z_\delta$'s with $k-1$ or fewer parts}  \\
&=& z_{g_i}z_{g_{i+1}\dots g_k g_1\dots g_{i-1}} + 
\textrm{a sum of $z_\delta$'s with $k-1$ or fewer parts}\\
&=& z_{g_i}z_{\lambda(g_{i+1}\dots g_k g_1 \dots g_{i-1})}+ 
\textrm{a sum of $z_\delta$'s with $k-1$ or fewer parts}\\
&=& z_{\lambda(\gamma)} +
\textrm{a sum of $z_\delta$'s with $k-1$ or fewer parts,}
\end{eqnarray*}
where the third equality uses the induction hypothesis. 
A trivial induction on the length of $\gamma$
now shows that any $z_{\gamma}, \gamma \vDash n$
can be reduced in the quotient to a linear combination of  $z_{\lambda}, 
\lambda
\vdash n$.

\end{proof} 
 
 
As a consequence of the proof we get that the 
ribbon Schur functions $r_{\lambda}$, $\lambda\vdash n$, 
span $\Lambda_{n}$ and so form a basis.  However, for general 
compositions $\alpha$ and $\beta$, it may be that $r_{\alpha}=r_{\beta}$. 
The rest of this section begins to deal with the question of when this 
can occur.  In principle, the relation $r_{\alpha}=r_{\beta}$ is a 
consequence of the relations \eqref{ribbonmult}, and more specifically 
\eqref{ribbonlinear}, so a purely algebraic development of the main results 
of this paper should be possible.  This has not yet been done.

\subsection{Equivalence of compositions} 
 
 We  define an algebraic 
equivalence on compositions and reinterpret it in a combinatorial manner. 
 
\begin{definition}\label{beta-equiv}Let $\beta$, $\gamma$ be 
compositions. We say $\beta$ and $\gamma$ are \emph{equivalent}, denoted 
$\beta \sim \gamma$, if for all 
$F=\sum c_\alpha F_\alpha \in \Lambda$, $c_\beta=c_\gamma$. 
\end{definition} 
\noindent 
That is, $\beta \sim \gamma$ if $F_\beta$ has the same coefficient 
as $F_\gamma$ in the expression of every symmetric function. 
Note that any basis for $\Lambda$ can be used as a finite test set 
for this equivalence. We will be particularly interested in the 
monomial symmetric function basis (\ref{monomial}) and the Schur 
function basis (\ref{schur}). 
 
\begin{example} For $\beta = 211$ and $\gamma = 121$ we find that 
$\beta \not\sim \gamma$ since 
$s_{22}= F_{22}+F_{121}$. 
\end{example} 
 
For any composition $\beta\vDash n$, we define $\M (\beta)$ to be the {\em multiset} of 
partitions determined by all coarsenings of $\beta$, that is, 
\begin{eqnarray} 
\M(\beta) = \{ \lambda(\alpha) ~|~ \alpha \ge \beta \}. 
\label{mbeta} 
\end{eqnarray} 
We denote by $\mult_{\M(\beta)}(\lambda)$ the multiplicity of $\lambda$ 
in $\M(\beta)$. 
 
\begin{example} 
Note that while 2111 and 1211 have identical {\em sets} of 
partitions arising from their coarsenings, $\M(2111)\not= \M(1211)$ 
since $\mult_{\M(2111)}(311)=1$ while $\mult_{\M(1211)}(311)=2$. 
\end{example} 
 
With this in mind we reformulate our equivalence. 
Recall that for $F\in\QQ$, 
$$F=\sum c_\beta F_\beta = \sum d_\beta M_\beta,$$ 
where the $c_\alpha$ and $d_\alpha$ are related by 
\begin{eqnarray} 
d _\beta =\sum  _{\alpha \ge \beta}c_\alpha,\quad 
   c_\beta =\sum _{\alpha\ge \beta} 
        (-1)^{l(\alpha)-l(\beta)}d_\alpha, 
\label{inversion} 
\end{eqnarray} 
and that $F\in\Lambda$ if and only if $d_\alpha = d_\beta$ whenever 
$\lambda(\alpha)=\lambda(\beta)$. 
The following is a 
direct consequence of (\ref{monomial}) and (\ref{inversion}). 
 
\begin{proposition} 
If the monomial symmetric function 
$m_\lambda = \sum_{\beta \models n} c_\beta F_\beta$, 
then 
$$c_\beta = [m_\lambda]_{F_\beta} = (-1)^{l(\lambda)-l(\beta)} 
 \mult_{\M(\beta)}(\lambda),$$ 
that is, up to sign, $c_\beta$ is the multiplicity 
of $\lambda$ in the multiset $\M(\beta)$. 
\label{monomulti} 
\end{proposition} 
 
As an immediate consequence we get 
 
\begin{corollary} 
\label{multi-schureq} 
If $\beta$ and $\gamma$ are 
compositions, then $\beta\sim\gamma$ if and only if 
$\M(\beta)=\M(\gamma)$. 
\end{corollary}

\begin{example} 
Returning to the example $\beta=211$ and $\gamma=121$, 
it is now straightforward to 
deduce $\beta\not\sim\gamma$ since 
\begin{eqnarray*} 
\M(\beta)=\{4, 31, 22, 211\}&\neq&\{4, 31, 31, 
211\}=\M(\gamma).\end{eqnarray*} 
\end{example} 
 
 
We are now ready to state the main result of this section. 
Note first that Proposition \ref{ribbonform} can be written 
\begin{equation}\label{ribbonform2} 
r_{\alpha}=  (-1)^{l(\alpha)} \sum_{\lambda \in \M(\alpha)}(-1)^{l(\lambda)} h_{\lambda}. 
\end{equation} 
 
\begin{theorem}\label{ribbon-equality} 
For the ribbon Schur functions $r_\beta$ and $r_\gamma$ corresponding 
to  compositions $\beta$ and $\gamma$, we have $r_\beta = r_\gamma$ 
if and only if $\M(\beta)=\M(\gamma)$. 
\end{theorem} 
\begin{proof} 
If $\M(\alpha)=\M(\beta)$, then by (\ref{ribbonform2}), $r_\beta = r_\gamma$. 
Conversely, since the $h_{k}$ are algebraically independent, equality of 
$r_\beta$ and $r_\gamma$ implies, again by  (\ref{ribbonform2}), that $\M(\alpha)=\M(\beta)$. 
\end{proof} 
 
An immediate corollary of this and (\ref{LRcoeff}) are the following Littlewood-Richardson 
coefficient identities. 
 
\begin{corollary} 
Suppose ribbon skew shapes $\lambda/\mu$ and $\rho/\eta$ correspond to 
compositions $\beta$ and $\gamma$, where $\beta \sim \gamma$ are both compositions of $n$. 
Then, for all partitions $\nu$ of $n$, 
$$c_{\mu,\nu}^\lambda = c_{\eta,\nu}^{\rho}.$$ 
\label{LRidents} 
\end{corollary}

\begin{example} 
Since $\M (211)=\{4, 31, 22, 211\}=\M (112)$ the above theorem 
assures us that 
$s_{222/11}=s_{4331/2221}$ 
and so, by Corollary \ref{LRidents}, 
$c^{222}_{11,\nu} = c^{4331}_{2221,\nu}$ 
for all partitions $\nu$ of 4. 
\end{example}

An immediate consequence of Theorem~\ref{ribbon-equality}, 
Corollary ~\ref{multi-schureq} and \cite[Corollary 7.23.4]{ECII} is another description 
of the equivalence $\sim$.  For $\sigma = \sigma(1)\sigma(2)\ldots \sigma(n) \in S_n$, 
the \emph{descent set} of $\sigma$ is defined to be the set 
$d(\sigma):=\{ i~|~ \sigma (i) > \sigma (i+1)\} \subset [n-1]$. 
 
\begin{corollary}\label{perm-descents} 
For $\beta, \gamma \vDash n$, $\beta \sim \gamma$ if and only if for all $\alpha \vDash n$, 
the number 
of permutations $\sigma\in S_n$ satisfying  $ d(\sigma)=S(\alpha)$ 
and $d(\sigma ^{-1})=S(\beta)$  is equal to the number of 
permutations $\sigma\in S_n$ satisfying  $ d(\sigma)=S(\alpha)$ 
and $d(\sigma ^{-1})=S(\gamma)$. 
\end{corollary} 
 
 
%
 

\section{Compositions of Compositions \label{composition-sec}} 
 
In this section we describe a method to combine compositions into larger ones 
that corresponds to determining the descent set of the tensor product of two permutations. 
This leads naturally to a necessary and sufficient condition for two compositions 
to be equivalent.

\subsection{Composition, tensor product, and
plethysm.} 
 
Let $\C _n$ denote the set of all compositions of $n$ and let 
$$ \C = \bigcup _{n\geq 1}\C _n.$$ 
Given $\alpha =\alpha _1 \ldots \alpha _k\vDash m$ and 
$\beta =\beta _1 \ldots \beta _l\vDash n$, recall the binary 
operations of \emph{concatenation} 
\begin{eqnarray*} 
\cdot : \C _m \times \C _n &\rightarrow&\C _{m+n}\\ 
(\alpha , \beta)&\mapsto& \alpha \cdot \beta = \alpha _1\ldots \alpha _k\beta _1\ldots 
\beta _l 
\end{eqnarray*} 
and \emph{near concatenation} 
\begin{eqnarray*} 
\odot : \C _m \times \C _n &\rightarrow&\C _{m+n}\\ 
(\alpha , \beta)&\mapsto& \alpha\odot\beta = \alpha _1\ldots \alpha _{k-1}(\alpha _k + 
\beta _1)\beta _2 \ldots \beta _l. 
\end{eqnarray*} 
For 
convenience we write 
\begin{equation*} 
\genfrac{}{}{0pt}{}{\alpha ^{\odot n} =}{} 
\genfrac{}{}{0pt}{}{\underbrace{\alpha \odot\alpha \odot\ldots 
\odot\alpha}}{n}. 
\end{equation*}

These two operations can be combined to produce a third, which 
will be our focus: 
\begin{eqnarray*} 
\circ : \C _m \times \C _n &\rightarrow&\C _{mn}\\ 
(\alpha , \beta)&\mapsto&\alpha\circ\beta = \beta ^{\odot\alpha _1}\cdot\beta ^{\odot\alpha 
_2}\cdots\beta ^{\odot\alpha _k}. 
\end{eqnarray*} 
 
\begin{example} If $\alpha = 12, \beta = 12$ then 
$\alpha\cdot\beta=1212$, $\alpha \odot \beta = 132$ and 
$\alpha\circ\beta = 12132$. 
\end{example} 
 
It is straightforward to observe that $\C$ is closed under $\circ$ 
and that for $\alpha \vDash m$ we have $1\circ \alpha = \alpha 
\circ 1 = \alpha$.  Note that the operation 
$\circ$ is not commutative since $12 \circ 3 = 36$ whereas 
$3\circ 12 = 1332$.

We now see that composing compositions corresponds to determining 
descent sets in the tensor product of permutations. 
 
\begin{definition}\label{tensor-prod} 
Let $\sigma =\sigma(1)\sigma(2)\ldots\sigma(m)\in S_m$ and 
$\tau =\tau(1)\tau(2)\ldots \tau(n)\in S_n$. Then their tensor product is 
the permutation 
\begin{eqnarray*} 
\sigma \otimes \tau &=& [(\sigma(1)-1)n+\tau(1)] 
[(\sigma(1)-1)n+\tau(2)]\ldots [(\sigma(1)-1)n+\tau(n)]\\ 
&&[(\sigma(2)-1)n+\tau(1)] \ldots [(\sigma(m)-1)n+\tau(n)]\quad 
\in S_{mn} .\end{eqnarray*} 
\end{definition} 
 
\begin{remark} An alternative realization is as follows. Given 
$\sigma \in S_m, \tau \in S_n$ and the $m\times n$ matrix 
$$M_{mn}=\begin{pmatrix} 
1&2&\cdots&n\\ 
n+1&n+2&\cdots&2n\\ 
\vdots&\vdots&\ddots&\vdots\\ 
(m-1)n +1&(m-1)n+2&\cdots&mn\end{pmatrix}$$ then $^\sigma M _{mn}$ 
is the matrix in which the $i$-th row of $^\sigma M _{mn}$ is the 
$\sigma (i)$-th row of $M_{mn}$. Similarly,  $M^\tau _{mn}$ 
is the matrix in which the $j$-th column of $M^\tau _{mn}$ is the 
$\tau (j)$-th column of $M_{mn}$. With this in mind, $\sigma 
\otimes \tau \in S_{mn}$ is the permutation obtained by reading 
the entries of $^\sigma M^\tau _{mn}$ by row. 
\end{remark} 
 
\begin{example} 
If $\sigma=213, \tau=132\in S_3$ then $\sigma\otimes\tau=213\otimes 132 = 465132798$, and 
$$M_{mn} =\begin{pmatrix} 
1&2&3\\ 
4&5&6\\ 
7&8&9\end{pmatrix}, ^\sigma M _{mn} =\begin{pmatrix} 
4&5&6\\ 
1&2&3\\ 
7&8&9\\ 
\end{pmatrix}, 
^\sigma M^\tau _{mn} =\begin{pmatrix} 
4&6&5\\ 
1&3&2\\ 
7&9&8 
\end{pmatrix}.$$\end{example}

The following shows that the operation $\circ$ on compositions yields 
the descent set of the tensor product of two permutations from their 
respective descent sets. 
 
\begin{proposition}\label{lemma-assoc} 
Let $\sigma \in S_m$ and $\tau \in S_n$. If $d(\sigma)=S(\beta)$ 
and $d(\tau)=S(\gamma)$ then $d(\sigma\otimes\tau)=S(\beta\circ 
\gamma)$. 
\end{proposition} 
 
\begin{proof} 
Let $d(\sigma)=S(\beta)= \{i_1,i_2,\ldots , i_k\}$ and 
$d(\tau)=S(\gamma)= \{j_1,j_2,\ldots , j_l\}$. Then 
$d(\sigma\otimes\tau)=\{j_1,j_2,\ldots j_l, n+j_1, n+j_2, \ldots , 
n+j_l,\ldots , (m-1)n+j_1,(m-1)n+j_2, \ldots , (m-1)n+j_l\} \cup 
\{ ni_1, ni_2, \ldots , ni_k\} = S(\beta\circ\gamma)$. 
\end{proof}

From Proposition \ref{lemma-assoc} and the associativity of $\otimes$, we can 
conclude that $\circ$ is associative.  Consequently we obtain 
 
\begin{proposition} 
$(\C , \circ )$ is a monoid. 
\end{proposition} 
 
\medskip

Finally we relate the operation $\circ$ on compositions to the operation of plethysm
on symmetric functions.
For the power sum symmetric function $p_m\in \Lambda_m$ and 
$g\in \Lambda_n$, define the {\em plethysm} 
$$p_m\circ g = p_m[g]=g(x_1^m,x_2^m,\dots).$$
Extend this to define $f\circ g\in \Lambda_{mn}$ for any 
$f\in \Lambda_m$ and $g\in \Lambda_n$ by requiring that
the map taking $f$ to $f\circ g$ be an algebra map. (See \cite[p.135]{MacD},
\cite[p.447]{ECII} for details.)  

When $\alpha \vDash m$ and $\beta \vDash n$ the ribbon functions
$r_{\alpha} \circ r_{\beta}$ and $r_{\alpha \circ \beta}$ both have degree $mn$.
While they are not in general equal, they are equal on the average, as seen
by the following identity.
                                                                              
\begin{proposition}                                                                              
For any {$\beta \vDash n$},
$${\sum_{\alpha \vDash m} r_{\alpha} \circ r_{\beta} =
 \sum_{\alpha \vDash m} r_{\alpha \circ \beta}} .$$
\end{proposition}

\begin{proof}
We first note that
\begin{equation*} r_{1}^{m } = \sum_{\alpha \vDash m} r_{\alpha}, \end{equation*}
which can be seen, for example, by repeated application of (\ref{ribbonmult}).  Since the plethysm
$f \circ g$ gives an algebra map in $f$, it follows that
$$\sum_{\alpha \vDash m} r_{\alpha} \circ r_{\beta} = r_{1}^{m}\circ r_{\beta} = r_{\beta}^{m}.$$

On the other hand, we show $r_{\beta}^{m} =  \sum_{\alpha \vDash m} r_{\alpha \circ \beta}$
by induction on $m$.  The case $m=1$ is clear.  For $m>1$, we get by induction
and (\ref{ribbonmult}) that
\begin{eqnarray*} 
 r_{\beta}^{m} &=& r_{\beta} \cdot  \sum_{\alpha \vDash m-1} r_{\alpha \circ \beta} \cr
 &=&  \sum_{\alpha \vDash m-1} \left[ r_{\beta \cdot (\alpha \circ \beta)} + 
 r_{\beta \odot (\alpha \circ \beta)}  \right] \cr
 &=&  \sum_{\alpha \vDash m} r_{\alpha \circ \beta}.
\end{eqnarray*} 
\end{proof}

\subsection{Unique factorization and other properties.} 
 
If a composition $\alpha$ is written in the form $\alpha _1\circ 
\alpha _2\circ \cdots \circ \alpha _k$ then we call this a 
\emph{decomposition} or \emph{factorization} of $\alpha$. 
A factorization $\alpha=\beta\circ\gamma$
is called \emph{trivial} if 
any of the following conditions are satisfied:
\begin{enumerate}
\item one of $\beta$, $\gamma$ is the composition 1,
\item the compositions $\beta$ and $\gamma$ both have length 1,
\item the compositions  $\beta$ and $\gamma$ both have all 
components equal to 1.
\end{enumerate}

\begin{definition} A factorization $\alpha=\alpha_1\circ\dots\circ \alpha_k$
is called \emph{irreducible} if no $\alpha_i\circ\alpha_{i+1}$ is a 
trivial factorization, and each $\alpha_i$ admits
only trivial factorizations.  In this case, each $\alpha_i$ is
called an \emph{irreducible factor}.
\end{definition}
 
\begin{theorem} \label{decomp-unique} 
The irreducible factorization of any composition is unique. 
\end{theorem} 
 
\begin{proof} 
We proceed by induction on the number of irreducible factors in a 
decomposition. 
 
First observe that  if the only irreducible factor of a 
composition is itself then its irreducible factorization is 
unique. 
 
Now let $\alpha$ be some composition with two irreducible 
factorizations 
$$\mu _1\circ\ldots\circ\mu_{k-1}\circ\mu _k = \alpha = \nu_1\circ\ldots\circ\nu_{l-1}\circ\nu 
_l,$$ and for convenience set $\beta = \mu 
_1\circ\ldots\circ\mu_{k-1}$, $\gamma = \mu _k$, 
$\delta =\nu_1\circ\ldots\circ\nu_{l-1}$ 
and $\epsilon = \nu _l$ so 
$$\beta\circ\gamma= \alpha =\delta \circ \epsilon.$$ 
Our first task is to establish $|\gamma |=|\epsilon |$ from which 
the induction will easily follow. First assume $|\gamma |=n$ and 
$|\epsilon |=s$ such that $s\neq n$ and without loss of generality 
let $s<n$. 
 
If $\epsilon =s$ then it follows $\gamma \neq n$ as if $\gamma =n$ then 
by our induction assumption and the fact that $s\neq n$ we 
have that the lowest common multiple of $s$ and $n$ would also be 
an irreducible factor, which is a contradiction. Hence $\gamma 
\neq n$ and so $l(\gamma )>1$. Furthermore since $l(\gamma)>1$ 
then $\gamma =\gamma _1\ldots \gamma _k$ must consist of 
components of $\alpha$ (the righthandmost and $k-1$ lefthandmost 
components, for example), which implies $s|\gamma _1,\ldots 
,s|\gamma _k$ and hence $\gamma $ is not an irreducible factor. 
 
Thus $\epsilon\neq s$ so $l(\epsilon )>1$ and since $s<n$ we also 
have that $l(\gamma )>1$. In addition,  since 
$l(\gamma)>1,l(\epsilon)>1$ we have as above that $\gamma$ and 
$\epsilon$ must consist of components of $\alpha$. Hence if $s|n$ 
then it follows that $\gamma$ has $\epsilon$ as an irreducible 
factor and hence $\gamma$ is not an irreducible factor. 
 
Consequently we have that if $s\neq n$ then $l(\gamma)>1, 
l(\epsilon)>1$ and $s\nmid n$. Moreover, the components of 
$\gamma$ consist of the components of $\epsilon$ repeated (and 
perhaps the sum of the first and last components of $\epsilon$) 
plus one copy of $\epsilon$ truncated at one end of $\gamma$. 
However, since $\gamma$ and $\epsilon$ consist of components of 
$\alpha$ it follows that if $s\neq n$, $l(\gamma)>1, 
l(\epsilon)>1$ and $s\nmid n$, then 
 $\epsilon$ cannot be an 
irreducible factor. Thus $|\gamma|=n=s=|\epsilon|.$ 
 
Now that we have established $|\gamma|=|\epsilon|$ we will show 
that in fact $\gamma = \epsilon$. If $\gamma =n$ then clearly 
$\epsilon =n$ and we are done. If not, then since the last 
component of $\beta,\delta \geq 1$ it follows the righthand 
components of $\alpha$ whose sum is less than $n$ must be those of 
$\gamma$ and $\epsilon$ and since $|\gamma|=|\epsilon|$ it follows 
that $\gamma=\epsilon$. 
 
Since we now have $\beta\circ\gamma = \alpha =\delta \circ 
\gamma$, it is straightforward to see $\beta =\delta$. By the 
associativity of $\circ$ the result now follows by induction. 
\end{proof}

We can also deduce expressions for the content and length of a 
composition in terms of its decomposition. We omit the proofs, 
which each follow by a straightforward induction.

\begin{proposition}\label{content-comp} 
For compositions $\beta _1,\beta _2,\ldots ,\beta_k$ 
$$|\beta _1\circ\beta _2\ldots\circ\beta_k|=\prod _{i=1}^k |\beta _i|.$$ 
\end{proposition}

\begin{proposition}\label{length-comp} 
For compositions $\beta _1,\beta _2,\ldots ,\beta_k$ 
$$l(\beta _1\circ\beta _2\ldots\circ\beta_k)=l(\beta _1)+\sum _{i=2}^k 
\left(\prod _{j=1}^{i-1} |\beta _j|\right)(l(\beta _i)-1).$$ 
\end{proposition}

Finally, it will be useful to observe that reversal of 
compositions commutes with the composition.  The proof is clear. 
 
\begin{proposition}\label{comp-reversal} 
Let $\beta$, $\gamma$ be compositions then 
$$(\beta\circ\gamma)^\ast = \beta ^\ast\circ\gamma ^\ast .$$ 
\end{proposition}

\begin{remark} 
For $\sigma \in S_n$, define $\sigma^\ast \in S_n$ by 
$\sigma ^\ast (i) := (n+1) -\sigma(n+1-i)$.  It is easy to see that 
for $\sigma\in S_m$, $\tau\in S_n$ 
$(\sigma\otimes\tau)^\ast = \sigma ^\ast\otimes\tau^\ast$ and 
$\beta(d(\sigma^\ast)) = \left(\beta(d(\sigma)\right)^\ast$. 
One wonders whether this, in conjunction with Proposition \ref{lemma-assoc}, 
can provide a more direct approach to that of the next two sections. 
\end{remark}

\section{Equivalence of Compositions under $\circ$ \label{decomp-proof}}

We show in this section that the equivalence relation of 
Definition \ref{beta-equiv} is related to the composition of 
compositions via reversal of terms.  In particular, we prove

\begin{theorem} 
\label{decomp-thm} Two compositions $\beta$ and $\gamma$ satisfy 
$\beta \sim \gamma$ if and only if for some $k$, 
$$\beta = \beta_1\circ\beta_2\circ\cdots\circ\beta_k 
\hbox{\rm \quad and \quad} 
 \gamma = \gamma_1\circ\gamma_2\circ\cdots\circ\gamma_k,$$ 
where, for each $i$, either $\gamma_i = \beta_i$ or $\gamma_i 
=\beta_i^\ast.$ Thus the equivalence class of a composition 
$\beta$ will contain $2^r$ elements, where $r$ is the number of 
nonsymmetric (under reversal) irreducible factors in the 
irreducible factorization of $\beta$. 
\end{theorem} 
 
Before we embark on the proof, which will consist of the remainder 
of this section and the next, we  note a corollary that follows 
immediately from Corollary~\ref{multi-schureq}, Theorem~\ref{ribbon-equality}, 
Corollary~\ref{perm-descents} and Theorem ~\ref{decomp-thm}.

\begin{corollary}\label{ribbonschur-equiv} 
The following are equivalent for a pair of compositions $\beta, \gamma$: 
\begin{enumerate} 
\item 
$r_\beta = r_\gamma$, 
\item 
in all symmetric functions $F=\sum c_\alpha F_\alpha$, the 
coefficient of $F_\beta$ is equal to the coefficient of 
$F_\gamma$, 
\item 
$\M(\beta)=\M(\gamma)$, 
\item 
the number of 
permutations $\sigma\in S_n$ satisfying  $ d(\sigma)=S(\alpha)$ 
and $d(\sigma ^{-1})=S(\beta)$  is equal to the number of 
permutations $\sigma\in S_n$ satisfying  $ d(\sigma)=S(\alpha)$ 
and $d(\sigma ^{-1})=S(\gamma)$  for all $\alpha$, 
\item 
for some $k$, 
$$\beta = \beta_1\circ\beta_2\circ\cdots\circ\beta_k 
\hbox{\rm \quad and \quad} 
 \gamma = \gamma_1\circ\gamma_2\circ\cdots\circ\gamma_k,$$ 
and, for each $i$, either $\gamma_i = \beta_i$ or $\gamma_i =\beta_i^\ast$. 
\end{enumerate} 
\end{corollary} 
 
\begin{example} 
Since $12132$ has irreducible factorization $12\circ 12$, Corollary 
~\ref{ribbonschur-equiv} assures us that 
$$r_{12132}=r_{13212}=r_{21231}=r_{23121}$$ 
and, moreover,  these are the only ribbon Schur functions equal to 
$r_{12132}$.  In addition, from 
$$s_{54221/311}=r_{12132} = r_{13212}=s_{54431/332},$$ 
we can conclude from (\ref{LRcoeff}) the identity of 
Littlewood-Richardson coefficients 
$$c_{311,\nu}^{54221} = c_{332,\nu}^{54431}$$ 
for all partitions $\nu$ of 9. 
\end{example} 
 
\subsection{Reversal implies equivalence.}\label{revimpequiv} 
We recall that for compositions $\beta$ and $\gamma$, 
$\beta\sim\gamma$ if and only if 
$\M(\beta)=\M(\gamma)$ by Corollary \ref{multi-schureq}. 
{}From this and Proposition \ref{comp-reversal} it is easy to conclude 
 
\begin{proposition} 
\label{first-equiv} For compositions $\beta$ and 
$\gamma_1,\dots,\gamma_k$, 
$$\beta^\ast \sim \beta$$ 
and 
$$\gamma_1\circ\gamma_2\cdots\circ\gamma_k \sim 
\gamma_1^\ast\circ\gamma_2^\ast\cdots\circ\gamma_k^\ast.$$ 
\end{proposition} 
 
We show now that reversal of any of the terms in a decomposition 
of $\beta$ yields a composition equivalent to $\beta$. 
 
\begin{theorem} 
\label{reverse-components} For any compositions $\beta$, $\gamma$ 
and $\alpha$, 
\begin{enumerate} 
\item 
$\beta^\ast\circ\gamma \sim \beta\circ\gamma$, 
\item 
$\beta\circ\gamma^\ast \sim \beta\circ\gamma$ and 
\item 
$\beta\circ\alpha^\ast\circ\gamma \sim 
\beta\circ\alpha\circ\gamma$. 
\end{enumerate} 
\end{theorem} 
\begin{proof} 
By definition, 
$$\beta\circ\gamma=\gamma ^{\odot\beta _1}\cdot\gamma ^{\odot\beta _2}\cdots\gamma ^{\odot\beta 
_k}$$ and 
$$\beta^\ast\circ\gamma=\gamma ^{\odot\beta _k}\cdots\gamma ^{\odot\beta _2}\cdot\gamma ^{\odot\beta 
_1}.$$ To prove (1), note that any coarsening $\delta$ of 
$\beta\circ\gamma$ that does not involve adding terms in different 
components $\gamma^{\odot\beta_i}$ clearly corresponds to a 
coarsening of $\beta^\ast\circ\gamma$ that has the same sorting 
$\lambda(\delta)$.  On the other hand, a coarsening that involves, 
say, combining terms in $\gamma^{\odot\beta_i}$ with terms of 
$\gamma^{\odot\beta_{i+1}}$ can be viewed as a coarsening of the 
first sort of 
$$(\beta_1,\dots,\beta_{i-1},\beta_i + \beta_{i+1}, \beta_{i+2},\dots, 
\beta_k)\circ \gamma,$$ which can be seen to correspond to one 
arising as a coarsening of $\beta^\ast\circ\gamma$. 
 
Assertions (2) and (3) follow from (1) and Proposition \ref{first-equiv} 
via 
$$\beta\circ\gamma^\ast \sim \beta^\ast\circ\gamma \sim \beta\circ\gamma$$ 
and 
$$\beta\circ\alpha^\ast\circ\gamma \sim \beta^\ast\circ\alpha\circ\gamma^\ast 
\sim \beta\circ\alpha^\ast\circ\gamma^\ast \sim 
\beta\circ\alpha\circ\gamma,$$ respectively. 
\end{proof} 
 
One direction in the assertion of Theorem \ref{decomp-thm} now 
follows from Theorem \ref{reverse-components}.  The remainder of 
this section and the next is devoted to the proof of the other 
direction. 
 
\subsection{Equivalence implies reversal.} 
 
In this subsection, we prove the converse to the result 
established in the previous subsection: namely, that if $\beta 
\sim \gamma$, then there is a factorization 
$\beta=\beta_1\circ\dots\circ\beta_k$ such that 
$\gamma=\gamma_1\circ\dots\circ\gamma_k$, where $\gamma_i=\beta_i$ 
or $\beta_i^*$. We achieve this via two theorems. The first of 
these is 
 
\begin{theorem}\label{hugh-2} 
Let $\beta\sim\gamma$, and $\beta=\delta\circ \epsilon$. Then 
$\gamma$ can be decomposed as $\zeta\circ\eta$ with 
$\zeta\sim\delta$ and $\eta\sim\epsilon$. \end{theorem} 
 
\begin{example} Let $\beta=13212$ and $\gamma=12132$. It is straightforward 
to check that these two compositions are equivalent. Note we have 
that $\beta=21\circ 12$. Theorem ~\ref{hugh-2} says that there 
should be a decomposition $\gamma=\zeta\circ\eta$ with $\zeta \sim 
21$ and $\eta \sim 12$. We observe that $\gamma=12\circ 12$ 
satisfies these conditions. 
\end{example} 
 
In order to prove Theorem ~\ref{hugh-2} we require two lemmas: 
 
\begin{lemma}\label{hughb1} Let $\beta =\delta\circ\epsilon$ where $\beta\vDash n$. 
Let $\epsilon$ have size $m$ and $p$ components. Let 
$\lambda=\lambda_1\dots\lambda_k$ be a partition of $n$ which 
occurs in $\mathcal{M}(\beta)$.  Let $\bar\lambda_i $ be the 
remainder when $\lambda_i$ is divided by $m$, and suppose that the 
sum of the $\bar\lambda_i$ is $m$. Then the number of non-zero 
$\bar\lambda_i$ is at most $p$. 
\end{lemma} 
 
\begin{proof} 
Reordering the parts of $\lambda$ if necessary, let 
$\lambda_1\dots \lambda_k$ be a composition of $n$ which is a 
coarsening of $\beta$.  Now consider the composition of $m$ given 
by $\bar\lambda _1\dots\bar\lambda_k $ (where we omit any zero 
components).  This composition is a coarsening of $\epsilon$, and 
thus has at most $p$ components.  \end{proof}

\begin{lemma}\label{hughb3} Let $\beta =\delta\circ\epsilon$ where 
$\beta\vDash n$ and $\epsilon \vDash m$, and let $\lambda$ be a 
partition of $m$, with $k$ parts.  Then 
$$\mult_{\mathcal{M}(\beta)}(\lambda,n-m)=(k-1+\mult_{\mathcal{M}(\beta)}(m,n-m)) 
\mult_{\mathcal{M}(\epsilon)}(\lambda).$$ \end{lemma} 
 
\begin{remark} 
Note that in the statement of the previous lemma and subsequently, 
when the context is unambiguous, we will refer to the multiplicity 
of a composition in the multiset of coarsenings of a composition 
when we intend the multiplicity of the partition determined by 
that composition. \end{remark} 
 
\begin{proof} Given a way to realize $\lambda$ from $\epsilon$, there are 
$k-1+\mult_{\mathcal{M}(\beta)}(m,n-m)$ corresponding ways to 
realize $(\lambda,n-m)$ from $\beta$: one must pick where to put 
in the $n-m$ component.  \end{proof} 
 
\noindent 
\textit{Proof of Theorem ~\ref{hugh-2}.} Let the size of 
$\epsilon$ be $m$. Write $q=n/m$. Let the number of components of 
$\epsilon$ be $p$. 
 
Define $\zeta\vDash q$ by setting $S(\zeta)=\{i\mid mi\in S(\gamma)\}$. 
Now 
$\mult_{\mathcal{M}(\zeta)}(\lambda)=\mult_{\mathcal{M}(\gamma)} 
m\lambda$ where we write $m\lambda$ for the partition obtained by 
multiplying all the parts of $\lambda$ by $m$. Similarly, 
$\mult_{\mathcal{M}(\delta)}(\lambda)=\mult_{\mathcal{M}(\beta)} 
m\lambda$. Thus, the equivalence of $\delta$ and $\zeta$ follows 
from that of $\beta$ and $\gamma$. 
 
Define $\eta_i \vDash m$, $i=0,\dots, q-1$, 
by setting 
$$S(\eta_i)=\{x\mid 0<x<m,\ x+im \in S(\gamma)\}.$$ 
We wish to show that all the $\eta_i$ are 
equal and equivalent to $\epsilon$. 
 
For any $0\leq i \leq q-1$, the number of components of $\eta_i$ 
is at most $p$: otherwise, consider the composition of $\gamma$ 
consisting of 
\begin{itemize} 
\item $im$ plus the first component of $\eta_i$, 
\item the remaining components of $\eta_i$ except the last, 
\item the last component of $\eta_i$ plus $(q-1-i)m$. 
\end{itemize} 
The partition corresponding to this composition appears in 
$\mathcal{M}(\gamma)$ but by Lemma ~\ref{hughb1}, it cannot appear 
in $\mathcal{M}(\beta)$, which is a contradiction. 
 
The cardinalities of $S(\beta)$ and $S(\gamma)$ must be the same, 
and we have already seen that $|S(\beta)\cap 
m\mathbb{Z}|=|S(\gamma)\cap m\mathbb{Z}|$. We know that 
$|S(\beta)\cap (\mathbb{Z}\setminus m\mathbb{Z})|=q(p-1)$, so the 
same must hold for $\gamma$.  Now, since each of the $\eta_i$ has 
at most $p$ components, each of the $\eta_i$ must have exactly $p$ 
components. 
 
We now need the following lemma: 
 
\begin{lemma}\label{hughb2} Let $\beta$, $\gamma$, and the 
$\eta_i$ be as already defined. Let $0\leq i < j \leq q-1$. Let 
$S(\eta_i)=\{a_1<\dots <a_{p-1}\}$ and $S(\eta_j)=\{b_1<\dots < 
b_{p-1}\}$.  Then $a_t \geq b_t$ for all $t$. 
\end{lemma} 
 
\begin{proof} If this were not so, let $\nu$ be the 
partition consisting of the following: 
\begin{itemize} 
\item $im$ plus the first component of $\eta_i$, 
\item the second through $t$-th components of $\eta_i$, 
\item $(j-i)m + b_t-a_t$, 
\item the $t+1$-th through $p-1$-th components of $\eta_j$, 
\item the last component of $\eta_j$ plus $(q-j-1)m$. 
\end{itemize} 
 
Now $\nu$ appears in $\mathcal{M}(\gamma)$ but by Lemma 
~\ref{hughb1} does not appear in $\mathcal{M}(\beta)$, a 
contradiction.  \end{proof}

Let $\mu$ be the partition of $m$ determined by $\epsilon$. Let $x 
= \mult_{\mathcal{M}(\beta)}(m,n-m)\in\{ 0,1,2\}$. The 
multiplicity of $(\mu,n-m)$ in $\mathcal{M}(\beta)$ is $p-1+x$. 
 
Now consider the possible occurrences of $(\mu,n-m)$ in 
$\mathcal{M}(\gamma)$.  If the $t$-th element of $S(\eta_0)$ 
coincides with the $t$-th element of $S(\eta_{q-1})$, then we have 
one possible occurrence of $(\mu,n-m)$ with $n-m$ as the $t+1$-th 
component.  Also, since by the equivalence of $\beta$ and 
$\gamma$, $x$ of $\{m,n-m\}$ are in $S(\gamma)$, there are $x$ 
possible occurrences of compositions realizing $(\mu,n-m)$ such 
that the $n-m$ part is either the first or the last component. 
However, there must be $p-1+x$ realizations of $(\mu,n-m)$, so all 
these possibilities must actually realize the partition. 
 
In particular, this shows that $S(\eta_0)$ and $S(\eta_{q-1})$ 
must coincide. Now, by Lemma ~\ref{hughb2}, all the 
$S(\eta_i)$ must coincide, and we can now denote all the $\eta_i$ 
by $\eta$. 
The equality of the $\eta_i$ (in particular, the equality of 
$\eta_0$ and $\eta_{q-1}$) means that we can apply the same 
argument as in Lemma ~\ref{hughb3} to show that for $\lambda$ a 
partition of $m$ with $k$ parts, 
$$\mult_{\mathcal{M}(\gamma)}(\lambda,n-m)=(k-1+x)\mult_{\mathcal{M}(\eta)}(\lambda).$$ 
 
The equivalence of $\beta$ and $\gamma$ also implies the 
multiplicities of $\lambda$ in $\mathcal{M}(\eta)$ and 
$\mathcal{M}(\epsilon)$ are equal for any $\lambda$ that is a 
partition of $m$, and hence that $\epsilon$ and $\eta$ are 
equivalent, as desired. This completes the proof of the theorem. 
\qed \medskip

The second theorem requires the concept of reconstructibility of a 
composition. 
 
\begin{definition} 
A composition $\beta$ is said to be \emph{reconstructible} if 
knowing $\mathcal{M} (\beta)$ allows us to determine $\beta$ up to 
reversal.  \end{definition} 
 
\begin{example} The composition $112$ is reconstructible 
because if $\beta$ is a composition satisfying 
$\mult_{\mathcal{M}(\beta)}\lambda(211)=1$ and 
$\mult_{\mathcal{M}(\beta)} \lambda(22)=1$, then $\beta=112$ or 
$\beta=211$.\end{example} 
 
\begin{theorem}\label{hugh-1} 
If $\beta \vDash n$ is not reconstructible, then $\beta$ 
decomposes as $\delta \circ \epsilon$, where neither $\delta$ nor 
$\epsilon$ have size 1.  \end{theorem} 
 
\noindent
\textit{Proof.} We establish this result 
by defining a function $h$ on $\beta$ and then proving that if 
$\beta$ is not reconstructible then $h$ is periodic with period 
$|\epsilon|>1$. This, in turn, yields our result. Since the proof 
of the periodicity of $h$ is somewhat technical we will state the 
pertinent lemmas but postpone their proofs until Section 
~\ref{technical}. Before we define $h$ we need a few other definitions. 
 
\begin{definition} With respect to a composition $\gamma \vDash n$, 
for any $0<i<n$, 
we say that $i$ is of \emph{type} 0, 1, or 2, depending on whether 
there are 0, 1, or 2 occurrences of the partition $(i,n-i)$ 
in $\mathcal{M}(\gamma)$ 
or, equivalently, if 0, 1, or 2 of 
$i,n-i$ are in $S(\gamma)$. For $i = n/2$, if $n/2$ is an integer, 
we say that its type is twice the number of occurrences of 
$(n/2,n/2)$ in $\mathcal{M}(\gamma)$.\end{definition} 
 
\begin{example} In the composition $11231$, 1, 4 and 7 are type 2, 
2 and 6 are type 1, 3 and 5 are type 0. 
\end{example} 
 
Fix  a composition $\beta$ of $n$. Let $A_i$ be the set 
of those elements of $[n-1]$ that are of type $i$ with respect to $\beta$. 
If $A_1=\emptyset$, then clearly $\beta$ is reconstructible.  Note 
that $A_1=\emptyset$ exactly when $\beta$ is symmetric under 
reversal. 
Now suppose $A_1\ne \emptyset$. Let $k$ be the least element of 
$A_1$. Reversing $\beta$ if necessary, we may assume that $k \in 
S(\beta)$, and $n-k\not\in S(\beta)$. 
 
\begin{definition} 
For $j\in A_1$, we say that $j$ is \emph{determined} if we can 
tell whether or not $j \in S(\beta)$ from $\mathcal{M}(\beta)$ and 
the knowledge that $k \in S(\beta)$. 
\end{definition} 
 
\begin{example} In the composition $12132$, $A_1=[8]$.  The determined 
elements are $\{1,2,4,5,7,8\}$.  3 and 6 are undetermined, because 
$12132 \sim 13212$, and $3 \in S(12132)$, while $3 \not \in 
S(13212)$, and the reverse is true of $6$.  \end{example}

\begin{definition} 
For $x,y \in [n-1]$, we say that they \emph{agree} if they are of 
the same type and either both or neither are in $S(\beta)$. 
\end{definition} 
 
\begin{remark} 
Note that this second condition follows from the first for $x,y$ 
of even type. 
\end{remark} 
 
We extend the notion of type to all $\mathbb{Z}$ 
by saying that multiples of $n$ are type 0, and otherwise, $x$ has the 
same type as $x \bmod n$.

If every element of $A_1$ is determined, then $\beta$ is 
reconstructible. Suppose $\beta$ is not reconstructible, so there 
are undetermined elements of $A_1$. Let us define $T_0$ to be the 
set of all integers that are undetermined, where we extend the 
notion of determinedness to all integers by saying that, in 
general, $x$ is determined if and only if $x\bmod n$ is 
determined. Let $t_0$ be the greatest common divisor of $T_0$. We 
are going to define inductively a collection $T_i$ of sets of 
integers. We will write $T_{\leq j}$ for the union of 
$T_0,\dots,T_j$.  Let $t_j$ be the greatest common divisor of 
$T_{\leq j}$. 
 
\begin{definition} 
For $i>0$, let $T_i$ be the set of $x$ not divisible by $t_{i-1}$, 
such that there is some $t \in T_{i-1}$ with $x$ and $x+t$ of even 
type and disagreeing. 
\end{definition} 
 
Clearly, only finitely many of the $T_i$ are non-empty. Let $s$ be 
the greatest common divisor of all the $T_i$.  By convention, set 
$t_{-1}=n$. 
 
We are now ready to define the function $h$ and state the results 
needed in order to analyze its periodicity. Let $g$ and $h$ be the 
functions defined on $\mathbb{Z}$ with respect to $\beta$ by 
 
$$h(x)=\left\{\begin{array}{rl} 0 & \text{ if $x$ is type 0}\\ 
                               1 & \text{ if $x$ is type 1 and $x\bmod n\in S(\beta)$}\\ 
                               -1& \text{ if $x$ is type 1 and $x\bmod n\not\in S(\beta)$}\\ 
                                2&\text{ if $x$ is type 2} 
\end{array}\right.$$ 
 
and 
 
$$g(x)=\left\{\begin{array}{rl} 0 & \text{ if $x$ is of even type}\\ 
                               1 & \text{ if $x$ is type 1 and $x\bmod n\in S(\beta)$}\\ 
                               -1& \text{ if $x$ is type 1 and $x\bmod n\not\in 
                               S(\beta)$}. 
\end{array}\right.$$ 
Consider the following three statements concerning 
the functions $g$ and $h$ and the sets $T_i$.

\begin{itemize} 
 
\item[$\Pb_i$:]  The function $g$ is $t_i$-periodic except at multiples of $t_i$. 
 
\item[$\Qb_i$:] The function $h$ is $t_{i-1}$-periodic except at multiples of 
$t_i$. 
 
\item[$\Rb_i$:] For $x\in T_{i+1}$ and $z$ of type 1, $t_{i} \nmid z$, 
$z$ and $x+z$ agree. 
 
\end{itemize} 
 
These statements are all defined for $i\geq 0$.  Note that $\Qb_0$ 
is immediate, by our conventional definition of $t_{-1}$. The 
remaining statements will follow by simultaneous induction. 
 
\begin{example}\label{runex} Consider the composition $\beta=132121332=213\circ12$. 
We can write out the values of $g$ and $h$ on $[18]$ as strings of 
18 characters, writing $+$ for 1, and $-$ for $-1$. 
\begin{eqnarray*}  g &= +-0+-++-0+--+-0+-0\\ h&=+-0+-++-2+--+-0+-0 \end{eqnarray*} 
Since $\beta \sim 312\circ12$, 6 and 12 are type 1 undetermined. 
In fact, $T_0\cap[18]=\{6,12\}$; $t_0=6$. We next observe that 3 
and 9 belong to $T_1$ because 3 and 3+6 (resp. 9 and 9+6) are of 
even type but disagree, and $6 \in T_0$.  In fact, 
$T_1\cap[18]=\{3,9\}$.  Hence $t_1=3$.  All the $T_i$ for $i>1$ 
are empty.  Thus $t_i=3$ for $i>1$. 
 
We now take a look at the meanings of $\Pb_i$ and $\Qb_i$ for this 
choice of $\beta$.  $\Pb_0$ says that $g$ is 6-periodic except at 
multiples of 6.  $\Qb_0$ says $h$ is 18-periodic except at multiples 
of 6.  $\Pb_1$ says that $g$ is 
3-periodic except at multiples of 3.  $\Qb_1$ says that $h$ is 
6-periodic except at multiples of 3.    $\Pb_2$ says nothing more 
than $\Pb_1$.  $\Qb_2$ says that $h$ is 3-periodic except at 
multiples of 3. 
\end{example} 
 
For clarity of exposition, we will divide the proof of the 
simultaneous induction 
 into several parts: 
 
\begin{itemize} 
 
\item Proof of $\Pb_0$ (Lemma ~\ref{P0}). 
\item Proof that $\Pb_j$ and $\Qb_j$ for $j\leq i$  imply $\Rb_i$ (Lemma ~\ref{PQR}). 
\item Proof that $\Rb_{i}$ and $\Pb_i$ imply $\Pb_{i+1}$ (Lemma ~\ref{RPP}). 
\item Proof that $\Pb_{i+1}$ and $\Qb_i$ imply $\Qb_{i+1}$ (Lemma ~\ref{PQQ}). 
 
\end{itemize} 
 
These four lemmas establish the simultaneous induction. 
 
Observe that for $i$ sufficiently large, $t_i=t_{i-1}=s$. Thus 
$\Qb_i$ implies that $h$ is $s$-periodic except at multiples of 
$s$. We now apply the following lemma: 
 
\begin{lemma}\label{newlemma} The composition $\beta \vDash n$ 
 has a decomposition 
$\beta=\delta\circ\epsilon$ with $|\epsilon|=p$ if and only if $p$ 
divides $n$ and the function $h$ determined by $\beta$ is 
$p$-periodic except at multiples of $p$.\end{lemma} 
 
\begin{proof} 
Suppose $\beta$ has such a decomposition.  It is clear that $p|n$. 
Write $h_\beta$ for the function determined by 
$\beta$, and $h_\epsilon$ for the function determined by $\epsilon$. 
For $x\in [n-1]$, 
$p\nmid x$, $h_\beta(x)=h_\epsilon(x \bmod p)$, 
which proves the desired 
periodicity. 
 
Conversely, suppose that $h_\beta$ has the desired periodicity. Define 
$\epsilon\vDash p$ by setting $h_\epsilon|_{[0,p-1]}=h_\beta|_{[0,p-1]}$. 
Define 
$\delta\vDash n/p$ by setting $h_\delta(x)=h_\beta(px)$.  It is then clear 
that $\beta=\delta\circ\epsilon$.   \end{proof} 
 
\begin{example} Continuing Example ~\ref{runex}, and applying Lemma 
~\ref{newlemma} to the assertion of $\Qb_2$, that $h$ is 3-periodic 
except at multiples of 3, we conclude that $132121332=\delta\circ\epsilon$ 
where $|\epsilon| =3$, which is indeed true, since $132121332=213\circ12$. 
\end{example} 
 
Returning to the proof of Theorem ~\ref{hugh-1}, we see that  an 
application of 
Lemma ~\ref{newlemma} implies that $\beta=\delta \circ \epsilon$, where 
$|\epsilon|=s$.  We have $s<n$ since $\beta$ is not reconstructible. 
Since also $s>1$ (see Lemma ~\ref{spos}), this factorization is 
non-trivial. This proves Theorem ~\ref{hugh-1}.\qed \medskip 
 
We are now in a position to prove our main result. 
\medskip 
 
\noindent 
\textit{Proof of Theorem ~\ref{decomp-thm}.}  If $\beta$ and 
$\gamma$ satisfy $\beta \sim \gamma$ then by Theorem 
~\ref{hugh-1}, we can factor $\beta=\beta_1\circ\dots\circ\beta_k$ 
where all the $\beta_i$ are reconstructible.  Applying Theorem 
~\ref{hugh-2} repeatedly, we find that 
$\gamma=\gamma_1\circ\dots\circ\gamma_k$, where 
$\gamma_i\sim\beta_i$.  However, since the $\beta_i$ are 
reconstructible, $\gamma_i\sim\beta_i$ implies that 
$\gamma_i=\beta_i$ or $\beta_i^*$. 
 
Conversely, if $\beta=\beta_1\circ\dots\circ\beta_k$ and 
$\gamma=\gamma_1\circ\dots\circ\gamma_k$ such that either $\gamma 
_i =\beta _i$ or $\gamma _i=\beta _i ^\ast$ then by Theorem 
~\ref{reverse-components} it follows that $\beta\sim\gamma$. 
 
Finally, observe that by Theorem ~\ref{decomp-unique} the 
equivalence class of $\beta$ contains $2^r$ elements where $r$ is 
the number of non-symmetric compositions under reversal in the 
irreducible factorization of $\beta$.\qed 
 
\section{Technical Lemmas\label{technical}} 
 
In this section we prove the technical lemmas which we deferred from the 
previous section. 
 We begin with a basic lemma which will be useful 
throughout this section. 
 
\begin{lemma}\label{cheating} 
Let $\beta\vDash n$, and let $\alpha=m\circ\beta$ for some $m>1$. 
Then: 
\begin{enumerate} 
 
\item $\mathcal{M}(\alpha)$ can be determined from $\mathcal{M}(\beta)$. 
 
\item If $n$ does not divide $x$, then $x$ has the same type with respect to 
$\alpha$ as $x \bmod n$ does with respect to $\beta$. 
 
\item The functions $g$ and $h$ determined by $\alpha$ and $\beta$ coincide. 
 
\item $t_i(\alpha)=t_i(\beta)$. 
 
\item Each of $\Pb_i$, $\Qb_i$ and $\Rb_i$ 
holds for $\alpha$ if and only if it holds for $\beta$. 
 
\end{enumerate}\end{lemma} 
 
\begin{proof} Suppose 
we know $\mathcal{M}(\beta)$.  We wish to determine 
$\mathcal{M}(\alpha)$. This is equivalent to determining the 
equivalence class of $\alpha$ with respect 
 to equivalence for compositions.  By Theorem 
~\ref {hugh-2}, the equivalence class of $\alpha$ consists exactly 
of those compositions which can be written as $m\circ\gamma$ with 
$\gamma\sim \beta$.  Thus, knowing $\mathcal{M}(\beta)$ suffices 
to determine $\mathcal{M}(\alpha)$. 
 
Observe that (2), (3), and (5) are immediate from the definitions. 
For (4), we have to verify that $x$ is determined for $\alpha$ if 
and only if $x \bmod n$ is determined for $\beta$.  Suppose $x 
\bmod n$ is determined for $\beta$.  That says exactly that all 
compositions in the equivalence class of $\beta$ agree at $x \bmod 
n$.  By Theorem ~\ref{hugh-2}, the equivalence class of $\alpha$ 
consists of the single-part partition $m$ composed with elements 
of the equivalence class of $\beta$, and therefore $x$ is 
determined for $\alpha$.  The converse follows the same way. 
\end{proof} 
 
The purpose of this lemma is that at any step in the simultaneous 
induction that proves $\Pb_i$, $\Qb_i$ and $\Rb_i$, we can replace 
$\beta$ by $m\circ\beta$ if we so desire.

\subsection{Proof of $\Pb_0$.} 
In this subsection we prove $\Pb_0$ (Lemma ~\ref{P0}).  We also prove 
Lemma ~\ref{t0ne1}, which will be necessary for our proof of 
Lemma ~\ref{spos}. 
 
Let the elements of $T_0\cap[n-1]$ be 
$m_1 < \cdots < m_l$. Let $r_i=\gcd(m_1,\dots,m_i)$.  Note that $m_1$ 
and $n-m_1$ are both in $T_0$, so $r_l$ divides $n$, and therefore 
$r_l$ coincides with $t_0$, the greatest common divisor of $T_0$. 
We begin with some lemmas. 
 
\begin{lemma}\label{lemma1} 
Suppose that $x$, $y$, and $x+y$ all lie in $A_1$.  Then from 
$\M(\beta)$ we can tell if $x$, $y$, and $n-(x+y)$ all agree, or 
if they don't all agree. \end{lemma} 
 
\begin{proof} If $x$, $y$, and $n-(x+y)$ agree, then $(x,y,n-(x+y))$ does 
not appear in $\M(\beta)$.  Otherwise, it does appear. 
\end{proof} 
 
\begin{lemma}\label{lemma2} Suppose $x$, $y$, and $x+y$ lie in $[n-1]$ and 
exactly two of them lie in $A_1$.  Then we can determine from 
$\M(\beta)$ whether or not they agree. \end{lemma} 
 
\begin{proof} The proof is similar to that of Lemma ~\ref{lemma1}, though 
there are more cases to check.  It is sufficient to check the 
cases: $x$ type 0 (and the others type 1); $x$ type 2; $x+y$ type 
0; $x+y$ type 2.  In each case, one sees that the multiplicity of 
$(x,y,n-(x+y))$ in $\M(\beta)$ depends on whether the two type 1 
points agree or disagree. \end{proof} 
 
\begin{definition} We say that a function $f$ defined on a set of integers 
including $[p-1]$ is 
\emph{antisymmetric} on $[p-1]$ 
if $f(x)=-f(p-x)$ for $0<x<p$.  \end{definition} 
 
\begin{lemma}\label{lemma3} Let $f$ be a function on $[d-1]$ 
which takes values $0$, 
$1$, $-1$, $*$, and suppose that there is some $c<d$, 
such that for all $x$ for which both sides are well-defined 
 
\begin{eqnarray} 
f(x) &=& -f(d-x)\label{e1}\\ 
f(x) &=& -f(c-x)\label{e2}\\ 
f(x) &=& f(c+x)\label{e3} 
\end{eqnarray} 
except that if either side equals $*$, the equation is not required to 
hold. 
Further, we require that the points where $f$ takes the value $*$ are either 
exactly the multiples of $c$ less than $d$, or else no points at all. 
Let $r$ be the greatest common divisor of $c$ and $d$.  Then on multiples 
of $r$, 
$f$ takes on only the values $0$ and (possibly) $*$.  On non-multiples 
of $r$, $f$ is $r$-periodic, and $f$ is antisymmetric on $[r-1]$. 
\end{lemma} 
 
\begin{proof} 
The proof is by induction. We first consider the base case, which 
is when $r=c$.  Periodicity is (\ref{e3}).  Antisymmetry is 
(\ref{e2}). Notice (\ref{e3}) also implies that $f$ is constant on 
multiples of $c$; by (\ref{e1}) this constant value is either $*$ 
or $0$. 
 
Now we prove the induction step. 
Let $(\ref{e1}')$, $(\ref{e2}')$, 
$(\ref{e3}')$ denote (\ref{e1}), (\ref{e2}), and (\ref{e3}), 
with $d$ replaced by $c$ and $c$ replaced 
by $d \bmod c$. 
It is easy to see that $(\ref{e1}')$, $(\ref{e2}')$, and $(\ref{e3}')$ 
follow from 
(\ref{e1}), (\ref{e2}), (\ref{e3}).  Also, $f$ restricted to $[c-1]$ never takes on the value 
$*$.  The desired results now follow by induction. 
\end{proof}

\begin{lemma}\label{period} For $1\leq i\leq l$, 
 
(i) $g$ is antisymmetric on  $[r_{i}-1]$, 
 
(ii) $g$ is $r_i$-periodic on $[m_{i}-1]$ except at multiples of $r_i$. 
\end{lemma} 
 
\begin{proof} The proof is by induction on $i$. 
We begin by proving the base case, which is when $i=1$. 
 
Suppose $0<x<m_1$.  By assumption, $x$ and $m_1-x$ are determined 
if they are type 1. Suppose one of them is of even type, and the 
other is type 1.  Then by Lemma ~\ref{lemma2}, we can determine 
$m_1$, contradiction.  Suppose that $x$ and $m_1-x$ are both type 
1. If they agree, Lemma ~\ref{lemma1} allows us to determine 
$m_1$, contradiction.  Hence they must disagree.  This establishes 
(i) in the base case. In the base case, (ii) is vacuous. 
 
Now we prove the induction step. 
For $i\geq 2$  define a function $g_i$ 
on $[m_{i}-1]$, as follows: 
$$ 
g_i(x)=\left\{ \begin{array}{ll} * & \text{ if $r_{i-1} | x$}\\ 
g(x) & \text{ otherwise.} 
\end{array}\right.$$ 
 
We wish to apply Lemma ~\ref{lemma3} to $g_i$, with $d=m_{i}$, 
$c=r_{i-1}$. 
If $0<x<m_{i}$, and neither $x$ nor $m_{i}-x$ is a multiple of 
$r_{i-1}$ (so in particular, neither is type 1 undetermined), then, as 
in the proof of the base case, $g_i(x)=-g_i(m_{i}-x)$.  This is condition (\ref{e1}). 
 
Suppose both $x$ and $m_{i-1}+x<n$ are type 1 and determined.  If 
they disagree (which means that $x$ and $n-(m_{i-1}+x)$ agree), 
then we can determine $m_{i-1}$, contradiction. Similarly, if one 
is of even type and the other is type 1 determined, we can 
determine $m_{i-1}$, again a contradiction. It follows that $g_i$ 
is $m_{i-1}$-periodic on $[m_i-1]$ except at multiples of 
$r_{i-1}$.  However, by induction, $g_i$ is $r_{i-1}$-periodic 
except at multiples of $r_{i-1}$ on $[m_{i-1}-1]$, so $g_i$ is 
$r_{i-1}$-periodic except at multiples of $r_{i-1}$ on $[m_i-1]$. 
This establishes condition (\ref{e3}). Condition (\ref{e2}) 
follows by the induction hypothesis. 
 
Thus, 
we can apply Lemma ~\ref{lemma3}. This proves the induction step, 
and hence the lemma. \end{proof}

\begin{lemma}\label{P0} $\Pb_0$ holds, that is to say, $g$ is 
$t_0$-periodic except possibly at multiples of $t_0$.  \end{lemma} 
 
\begin{proof} 
Since $t_0=r_l$, we have already shown (Lemma ~\ref{period}) 
that $g$ is $t_0$-periodic on 
$[m_l-1]$ except at multiples of $t_0$. 
Since $g$ is antisymmetric on $[n-1]$ (by the definition of 
$g$) it follows that $g$ is $t_0$-periodic on $[n-1]$ except at multiples 
of $t_0$, from which the desired result follows. 
\end{proof} 
 
We now prove Lemma ~\ref{t0ne1} which will be used in the proof of 
Lemma ~\ref{spos}. 
 
\begin{lemma}\label{t0ne1} The greatest common divisor $t_0$ of $T_0$ does 
not divide $k$.  \end{lemma} 
 
\begin{proof} 
Suppose otherwise. Let $i$ be the least index such that $r_i | k$. 
Note that $i>1$, since $k<m_1$. By the result of applying Lemma 
~\ref{lemma3} to $g_i$, we know that $g_i$ is zero on multiples of 
$r_{i}$ which are not multiples of $r_{i-1}$.  However, this means 
that $g_i(k)=0$, which contradicts the fact that $k$ is type 1. 
\end{proof}

\subsection{Proof of $\Rb_i$.} 
We begin by deducing $\Rb_0$ from $\Pb_0$ (Lemma ~\ref{R0}).  We then 
prove the general statement that $\Pb_j$ and $\Qb_j$ for $j \le i$ 
 imply $\Rb_i$ (Lemma ~\ref{PQR}), 
which reduces to the argument for Lemma ~\ref{R0}. 
 
\begin{lemma}\label{R0} $\Pb_0$ implies $\Rb_0$.  \end{lemma} 
 
\begin{proof} We must show that if $z$ is type 1 and not a multiple of 
$t_0$ (which means in particular that it is determined), and $y_1 \in 
T_1$, then $z$ and $z+y_1$ agree.

Since $y_1\in T_1$, there is some $y_0\in T_0$ such that $y_1$ and 
$y_1+y_0$ are of even type and disagree. Clearly, we may assume that 
$z$, $y_0$, and $y_1$ are all positive. 
 
If $z+y_1+y_0>n$, we may replace $\beta$ by $m\circ \beta$ for 
some sufficiently large $m$, by Lemma ~\ref{cheating}.  We also 
wish to assume that $z<y_0$.  If this is not true, we can make it 
true by another replacement as above, followed by adding $n$ to 
$y_0$. 
 
By $\Pb_0$, we know that $z$ and $z+y_0$ agree.  Also, observe 
that since $z$ is type 1, so is $n-z$, and thus, by $\Pb_0$, so is 
any $w\equiv -z \bmod {t_0}$.  Since $y_1$ is of even type, this 
means that $t_0 \nmid z+y_1$, so $z+y_1$ is of even type or 
determined, and $\Pb_0$ tells us that $g(z+y_1+y_0)=g(z+y_1)$. 
 
By considering the multiplicity of $(y_0,y_1,z,n-(y_1+y_0+z))$ in 
$\mathcal{M}(\beta)$, we see that one of  two things happens: 
\begin{itemize} 
\item $z+y_1$ and $z+y_1+y_0$ are type 1 and both agree with $z$ 
\item $z+y_1$ agrees with $y_1$, while $z+y_1+y_0$ agrees with $y_0+y_1$. 
\end{itemize} 
 
We now exclude the second possibility.  Suppose we are in that 
case. Let $w=y_0-z$.  This $w$ is not a multiple of $t_0$, so $w$ 
is of even type or is determined. As already remarked, since $w \equiv 
-z \bmod t_0$, $w$ must be type 1 determined. 
Now apply the previous part of the 
proof to $(y_0',y_1',z')$ with $y_0'=y_0$, $y_1'=z+y_1$, $z'=w$. 
Then we see that either $w+z+y_1$ must either be the same type as 
$y_1$, or as $w$.  However, $w+z+y_1=y_0+y_1$, and we know that it 
is of even type but disagrees with $y_1$, which is a 
contradiction. 
\end{proof} 
 
\begin{lemma}\label{PQR} $\Pb_{j}$ and $\Qb_{j}$ for $j\leq i$ 
imply $\Rb_{i}$. 
\end{lemma} 
 
\begin{proof} 
We wish to show that for $z$ of type 1, $t_{i}\nmid z$ (so in particular 
$z$ is determined), and $x\in T_{i+1}$, that $z$ and $z+x$ agree. 
Write $y_{i+1}$ for $x$. 
Now there is some $y_{i}\in T_{i}$ such that $y_{i+1}$ and $y_i+y_{i+1}$ 
are of even type and disagree.  Similarly, choose 
$y_j$ for all $0\leq j 
\leq i-1$ so that $y_{j+1}$ and $y_{j+1}+y_j$ are even type and disagree.

For $I$ a subset of $[0,i+1]$, write $y_I$ for the sum of the 
$y_j$ with $j\in I$.  We now determine the types of $y_I$ and $y_I+z$. 
 
\begin{lemma}\label{types} Let $I \subset [0,i+1]$. Let $j$ be the maximal 
element of $I$.  Then: 
\begin{enumerate} 
\item If $I$ does not contain $j-1$, then $y_I$ agrees with $y_j$. 
\item If $I$ does contain $j-1$, then $y_I$ is of even type disagreeing with 
$y_j$. 
\item Either $z+y_I$ is of even type or it is determined. 
\item If $I$ does not contain $i+1$, $z+y_I$ agrees with $z$. 
\item If $I$ contains $i+1$ but not $i$, $z+y_I$ agrees with $z+y_{i+1}$. 
\item If $I$ contains $i+1$ and $i$, $z+y_I$ agrees with $z+y_{i+1}+y_{i}$. 
\item Either $z+y_{i+1}$ and $z+y_{i+1}+y_{i}$ agree, or they are both of even type. 
\end{enumerate} \end{lemma} 
 
\begin{proof} 
Statement (1) follows from $\Qb_{j-1}$, since $y_j$ is not a multiple of 
$t_{j-1}$. 
Statement (2) follows because $y_j$ and $y_j+y_{j-1}$ are 
of even type and disagree, and then applying $\Qb_{j-1}$ as before.

Since $g$ is $t_i$-periodic except at multiples of $t_i$, and its 
period is anti-symmetric, it follows that any $w\equiv -z \bmod 
{t_i}$ must be of odd type.  Thus $y_{i+1}\not\equiv -z \bmod 
t_i$.  All the other $y_l$ are multiples of $t_i$.  Thus $t_i\nmid 
z+y_I$, so $z+y_I$ is either determined or of even type.  This 
establishes (3). 
 
Statement (4) follows from $\Pb_{i}$, since $z$ is not a multiple of 
$t_{i}$. 
Since $t_i\nmid z+y_{i+1}$, (5) follows from $\Qb_i$. 
Statement (6) follows from $\Qb_{i}$ together with the fact that, since 
$t_{i}$ does not divide $z+y_{i+1}$, it doesn't divide $z+y_{i+1}+y_{i}$. 
Statement (7) follows from $\Pb_{i}$. 
\end{proof} 
 
We now return to the proof of $\Rb_i$. We want to assume that 
$p=n-(z+\sum_{j=0}^{i+1}y_{j})>0$, and that $p$ does not coincide 
with any $y_j$ or $z$.  In order to guarantee this, by 
Lemma~\ref{cheating}, we may replace $\beta$ by $m\circ\beta$, and 
add multiples of $n$ as desired to the $y_i$ and $z$. 
 
Since $y_0\in T_0$, it is undetermined.  This means precisely that there is 
some composition $\gamma$ which is equivalent to 
$\beta$ (but not equal to $\beta$), such that 
$k\in S(\gamma)$, but $y_0$ is in exactly one of $S(\beta)$, $S(\gamma)$. 
Note that since $\gamma$ is equivalent to $\beta$, every $0<x<n$ has the 
same type in $\beta$ and $\gamma$.

Write $\nu$ for the partition of $n$ 
whose parts are $(y_0,y_1,\dots, y_{i+1},z,p)$. 
One consequence of the equivalence of 
$\beta$ and $\gamma$ that we shall focus on is the fact that 
$\mult_{\mathcal{M}(\beta)}(\nu)=\mult_{\mathcal{M}(\gamma)}(\nu)$. 
 
Let $\Omega$ be the set of 
all the compositions of $n$ determining the partition 
$\nu$.  It will 
be convenient for us to keep track of such a composition as two lists: the 
left list, which consists 
of the components in order which precede $p$, 
and the right list, which consists 
of the components following $p$ in reverse order.  For any composition in 
$\Omega$, each 
component other than $p$ occurs in exactly one list, 
and any pair of lists with 
this property determines a composition. 
 
We put an order $\prec$ on the components 
$y_j$, $z$ by ordering the $y_j$ by their 
indices, and setting $y_j\prec z$ for $j \ne {i+1}$.  (Thus, 
the order is nearly 
a total order but not quite: $y_{i+1}$ and $z$ are incomparable.) 
 
\begin{definition} A composition in $\Omega$ is called \emph{ordered} if both its 
right and left lists are in (a linear extension of) $\prec$ order. 
The other compositions in $\Omega$ are called \emph{disordered}.\end{definition} 
 
\begin{lemma}\label{disordered} 
The number of disordered compositions which can be 
obtained as coarsenings of $\beta$ is the same as the number that can be 
obtained as coarsenings of $\gamma$.  \end{lemma} 
 
\begin{proof} 
To prove this lemma, we 
will define an involution $i$ on disordered compositions 
such that $\kappa$ is a coarsening of $\beta$ if and only if $i(\kappa)$ is a 
coarsening of $\gamma$. 
 
Fix a disordered composition $\kappa$. Let $M(\kappa)$ be the 
maximal subset of $y_0,\dots, y_{i+1}, z$ which is a $\prec$ order 
ideal such that $M(\kappa)$ consists of the union of initial 
subsequences of the left and right lists of $\gamma$, and these 
subsequences are in $\prec$ order. Write $M_L(\kappa)$ and 
$M_R(\kappa)$ for these two initial subsequences. Then $i(\kappa)$ 
is obtained by swapping $M_L(\kappa)$ and $M_R(\kappa)$. Observe 
that $i(\kappa)$ is disordered if and only if $\kappa$ is 
disordered. 
 
\begin{example} We give an example of the definition of $i$. 
$$ \mbox{If }\kappa=\begin{array}{|c|c|} 
y_0 & y_1\\ 
y_2& y_4\\ 
y_5 & y_3 \\ 
z& \end{array} \mbox{ then $M(\kappa)=\{y_0,y_1,y_2\}$ and } 
i(\kappa)=\begin{array}{|c|c|} 
y_1&y_0\\ 
y_5&y_2\\ 
z&y_4\\ 
&y_3 \end{array}.$$ \end{example} 
 
We shall now define a bijection, also denoted $i$, taking 
$S(\kappa)$ to $S(i(\kappa))$, such that for $x\in S(\kappa)$, 
$x\in S(\beta)$ if and only if $i(x)\in S(\gamma)$.  The existence 
of such a bijection between $S(\kappa)$ and $S(i(\kappa))$ implies 
that $\kappa$ is a coarsening of $\beta$ if and only if 
$i(\kappa)$ is a coarsening of $\gamma$, proving the lemma. 
 
To define the bijection between $S(\kappa)$ and $S(i(\kappa))$, 
we need another definition: 
 
\begin{definition}  We say $x\in S(\kappa)$ is an \emph{outside break} if it 
is either the sum of an initial subsequence of $M_L(\kappa)$ or $n$ minus 
the sum of an initial subsequence of $M_R(\kappa)$.  Otherwise, we say that 
$x\in S(\kappa)$ is an \emph{inside break}. \end{definition} 
 
\begin{example} In our continuing example, 
the outside breaks of $\kappa$ are $y_0$, $y_0+y_2$, and 
$n-y_1$, while the outside breaks of $i(\kappa)$ are 
$n-y_0$, $n-(y_0+y_2)$, and $y_1$. 
The inside 
breaks in $\kappa$ are $y_0+y_2+y_5$, $y_0+y_2+y_5+z$, $n-(y_1+y_4)$, 
$n-(y_1+y_4+y_3)$, while the corresponding inside breaks in $i(\kappa)$ are 
$y_1+y_5$, $y_1+y_5+z$, $n-(y_0+y_2+y_4)$, $n-(y_0+y_2+y_4+y_3)$. 
\end{example} 
 
If $x$ is an outside break of $\kappa$, set $i(x)=n-x$. Clearly, 
$i(x)$ is an outside break of $i(\kappa)$.  Now observe that all 
the outside breaks except $y_0$ or $n-y_0$ are of even type in 
$\beta$ by Lemma ~\ref{types}.  Thus for these outside breaks 
(excluding $y_0$ and $n-y_0$), $x\in S(\beta)$ if and only if $x$ is type 2 
for $\beta$ if and only if $x$ is type 2 for $\gamma$ if and only if $i(x)$ is type 2 
for $\gamma$ if and only if $i(x) \in S(\gamma)$.  On the other hand, $y_0 
\in S(\beta)$ if and only if $n-y_0 \in S(\gamma)$. Thus, for $x$ an outside 
break of $\kappa$, $x\in S(\beta)$ if and only if $i(x)\in S(\gamma)$.

Now we consider the inside breaks. 
Let $y^L$ denote the sum of the $y_j$ appearing in $M_L(\kappa)$, 
and similarly for $y^R$. 
If $x$ is an inside break for $\kappa$, 
set $i(x)=x-y^L+y^R$. This is clearly an inside break for 
$i(\kappa)$. 
 
Since $\kappa$ is disordered, define $l$ by $M(\kappa)=\{y_0, 
y_1,\dots,y_l\}$. By definition, all the $y_j$ that occur in $y^L$ 
and $y^R$ have $j\leq l$. To show that $x\in S(\beta)$ if and only 
if $i(x) \in S(\gamma)$ there are a four cases to consider: when 
$x$ is of the form $y_I$, $z+y_I$, $n-y_I$, or $n-y_I-z$.  In the 
first case, observe that $I$ contains at least one element greater 
than $l+1$, and so, by Lemma ~\ref{types}(1) or (2), $x$ and 
$i(x)$ agree and are of even type. It follows that $x\in S(\beta)$ 
if and only if $i(x)\in S(\beta)$ if and only if $i(x)\in 
S(\gamma)$, as desired. 
 
In the second case, since $l\leq i-1$ 
it is again clear  by Lemma 
~\ref{types}(4), (5), or (6), that $x$ and $i(x)$ agree, so $x\in 
S(\beta)$ if and only if $i(x)\in S(\beta)$. By Lemma 
~\ref{types}(3), $i(x)$ is either determined or of even type, so 
$i(x)\in S(\beta)$ if and only if $i(x)\in S(\gamma)$, which 
establishes the desired result. 
 
The third and fourth cases are similar to the first and second 
cases. This completes the proof that $i$ is a bijection from 
$S(\beta)$ to $S(\gamma)$, which completes the proof of the lemma. 
\end{proof}

Now we consider the ordered compositions.  Suppose $\kappa$ is an 
ordered composition which is a coarsening of $\beta$. Thus $y_0$ 
is the beginning of one list.  Which list is determined by which 
of $y_0$ and $n-y_0$ is a break in $\beta$. Since $y_1$ and 
$y_1+y_0$ disagree, which list $y_1$ occurs in is forced. 
Similarly for $y_2$, etc.  Hence all the $y_j$ are forced up to 
$y_{i}$.  There are now six possible ways to complete the 
construction.  For each of these six possibilities we show the 
positions of $y_{i}$, $y_{i+1}$, and $z$ in the two lists. 
$$\begin{array}{||l|l||l|l||l|l||l|l||l|l||l|l||} 
\multicolumn{2}{||c||}{(a)}& \multicolumn{2}{c||}{(b)}& \multicolumn{2}{c||}{(c)}& 
\multicolumn{2}{c||}{(d)} & \multicolumn{2}{c||}{(e)} & \multicolumn{2}{c||}{(f)} \\ 
\hline 
y_{i}& \qquad &y_{i}&z       & y_{i}& \qquad&y_{i}&z       &y_{i}&y_{i+1}       &y_{i}&y_{i+1} \\ 
z      &        & \qquad&y_{i+1}      &y_{i+1}     &       &y_{i+1}    &\qquad &\qquad &z         &z        &\qquad\\ 
y_{i+1}    &        &       &\qquad   & z      &       &\qquad &        &       &\qquad   &\qquad   & 
\end{array}$$ 
 
The argument now proceeds as in Lemma ~\ref{R0}. 
Essentially what has happened is 
that by reducing to ordered compositions, we do not need to consider the 
$y_j$ with $j<{i}$.  We are now only interested in the 
middle part of the composition, which involves parts $y_{i}$, $y_{i+1}$, 
$z$, and $p$.  Also $y_{i}$ now behaves like $y_0$ in Lemma ~\ref{R0}: 
we count 
up the number of compositions which occur with $y_{i}$ on the extreme 
left (among 
the four parts we are interested in) and those where it occurs on the extreme 
right.  One of these numbers represents the contribution of ordered partitions 
to $\mult_{\mathcal{M}(\beta)}(\nu)$, the other the contribution to 
$\mult_{\mathcal{M}(\gamma)}(\nu)$.  These numbers must therefore be the 
same. 
  As in the proof of Lemma ~\ref{R0}, we consider cases based on 
the types (and for type 1, whether or not each is a break) of 
$y_{i+1}$, $z$, $y_{i+1}+z$, and $y_{i+1}+y_{i}+z$. Lemma 
~\ref{types}(7) eliminates a number of possibilities and with the 
remainder, as in Lemma ~\ref{R0}, one of the following two things 
must happen: 
\begin{itemize} 
\item $z+y_{i+1}$ and $z+y_{i+1}+y_i$ are type 1 and both agree with $z$ 
\item $z+y_{i+1}$ agrees with $y_{i+1}$, while $z+y_{i+1}+y_{i}$ agrees with $y_{i+1}+y_{i}$. 
\end{itemize} 
We now exclude the second possibility. 
 
Since $z$ is not a multiple of $t_i$, $\Pb_i$ tells us that 
$y_i-z$ agrees with $n-z$, which is type 1 determined. Also by 
$\Pb_i$, $z\not\equiv-y_{i+1}\pmod {t_i}$, so $t_i$ does not 
divide $z+y_{i+1}$. Since $y_i\in T_i$, and $z+y_{i+1}$ and  
$z+y_{i+1}+y_i$ disagree, $z+y_{i+1}\in T_{i+1}$.   
Set $z'=y_i-z$, $y_{i+1}'=z+y_{i+1}$, 
$y_j'=y_j$ for $j\leq i$.  Applying the whole proof of the lemma 
so far, we find that $z'+y_{i+1}'$ must agree either with $z'$ or 
$y_{i+1}'$, which is to say that $y_i+y_{i+1}$ agrees with either 
$y_i-z$ or $z+y_{i+1}$, both of which are impossible, and we are 
done. 
\end{proof}

\subsection{Proofs of $\Pb_i$ and $\Qb_i$.} 
 
We begin with a preliminary lemma which will be useful for the proofs 
of Lemmas ~\ref{RPP} and ~\ref{PQQ}. 
While working towards proving these two lemmas, we will often need to consider 
$\mathbb{Z}/t_j\mathbb{Z}$ (for some $j$).  We will write $\mathbb{Z}_{t_j}$ 
for $\mathbb{Z}/t_j\mathbb{Z}$, and 
$\bar z$ for the image 
of $z$ in $\mathbb{Z}_{t_j}$. 
 
\begin{lemma}\label{prelim} 
Let $f$ be a function defined on $\mathbb{Z}$. Let $S\subset \ztj$ 
be such that $f(z)=f(z+t_{j-1})$ 
for $\bar z\in S$. 
Suppose further that for any 
$x \in T_j$, $f(z)=f(z+x)$ provided $\bar z\in S$. 
Then $f(z)=f(z+t_j)$ for all $\bar z\in S$. 
\end{lemma} 
 
\begin{proof} Write $t_j$ as the sum of a series 
of elements of $T_{\leq j}$.  Let 
the partial sums of this series be $x_1,x_2,\dots,x_m=t_j$.  Then 
observe that if $\bar z\in S$, then the same is true for 
$\overline{z+x_l}$ for all $l$.  It follows from the assumptions of the lemma that 
$f(z+x_l)=f(z+x_{l+1})$, and the result is proven. 
\end{proof} 
 
\begin{lemma}\label{neg} For $p>0$, if $x\in T_{p}$, then there is an element $x'\in T_p$ 
such that $x \equiv -x'$ modulo $t_{p-1}$.  \end{lemma} 
 
\begin{proof} Since $x \in T_p$, there is some $y\in T_{p-1}$ such that 
$x$ and $y+x$ are of even type and disagree.  It follows that $n-y-x$ and 
$n-x$ are of even type and disagree, and hence that $n-y-x \in T_p$. 
Set $x'=n-y-x$.  \end{proof}

\begin{lemma}\label{RPP} 
$\Rb_i$ and $\Pb_i$ imply $\Pb_{i+1}$.  \end{lemma} 
 
\begin{proof} Let $j=i+1$. 
We wish to show that $g$ is $t_{j}$-periodic except at 
multiples of $t_{j}$.  Let $S=\ztj 
\setminus \{\bar 0\}$.  $\Pb_i$ tells us that $g$ is $t_{j-1}$-periodic 
except at multiples of $t_{j-1}$. 
Suppose $\bar z \in S$, and $x \in T_j$. 
$\Rb_i$ tells us that if $z$ is type 1, then $g(z+x)=g(z)$. 
Likewise, if $z+x$ is type 1, then, choosing $x'$ as provided 
by Lemma~\ref{neg}, $z+x+x'$ is type 1, and 
now by the $t_{j-1}$ periodicity of $g$, $g(z+x)=g(z)$.  If neither 
$z$ nor $z+x$ is type 1, then $g(z+x)=0=g(z)$.

Thus, it follows that for any $z$ such 
that $\bar z\in S$, and $x$ in $T_j$, that $g(z+x)=g(z)$. 
Therefore, we can apply Lemma ~\ref{prelim}, 
and desired result follows.  \end{proof} 
 
\begin{lemma}\label{PQQ} 
$\Pb_{i+1}$ and $\Qb_{i}$ imply $\Qb_{i+1}$.  \end{lemma} 
 
\begin{proof} Let $j=i$.  Let $S=\ztj\setminus t_{j+1}\ztj$. 
We wish to show that $h(z)=h(z+t_j)$ for $\bar z \in S$. 
$\Qb_{i}$ tells us that $h(z+t_{j-1})=h(z)$ for $\bar z\in S$. 
Now suppose that we have some $z$ such that $\bar z \in S$, and 
$x \in T_{j}$.  By $\Pb_{i+1}$, if $h(z)= \pm 1$ then 
$h(z+x)=h(z)$.  Also by $\Pb_{i+1}$, if $h(z)$ is even, then so is 
$h(z+x)$. 
Now, if $h(z)\ne h(z+x)$, then 
$z\in T_{j+1}$, contradicting our assumption. Thus $h(z+x)=h(z)$ and we can 
apply Lemma ~\ref{prelim} to obtain the desired result. 
\end{proof} 
 
\subsection{Proof that $s>1$.}  Finally, we show that $s$, the greatest common 
divisor of the $T_i$, is greater than 1. 
 
\begin{lemma}\label{abel} 
Let $G$ be an arbitrary finite abelian group, which we write 
additively.  Let $Y$ be a 
set of generators for $G$, closed under negation.  Fix some $a\in G$. 
For any $b$ in $G$, it is possible to write $b$ as the sum of a 
series of elements from $Y$, so that no proper 
partial sum of the series equals 
$a$ (i.e., excluding the empty partial sum and the complete partial sum). 
\end{lemma} 
 
\begin{proof} The proof is by induction on $|G|$.  If $G$ is cyclic, pick 
$x\in Y$ a generator for $G$.  If $b$ occurs before $a$ in the sequence 
$x, 2x, \dots$, then we are done.  Otherwise, use $-x$. 
 
If $G$ is not cyclic, find a cyclic subgroup $H$ which is a direct 
summand, and has a generator $x\in Y$.  Let $\bar a$, $\bar b$ 
denote the images of $a$ and $b$ in $G/H$. Apply the induction 
hypothesis to $G/H$.  Lifting to $G$, we obtain a series whose sum 
differs from $b$ by an element of $H$, which we can dispose of as 
in the cyclic case above. The only problem occurs if $\bar b=\bar 
a, b\ne a$, and the series for $G/H$ happens to sum to $a$.  In 
this case, instead of putting the series obtained for $H$ after 
the series for $G/H$, begin with the first term from the series 
for $H$, followed by the series for $G/H$, followed by the rest of 
the series for $H$. 
\end{proof} 
 
\begin{example} Let $G\cong\mathbb{Z}/2\mathbb{Z}\oplus \mathbb{Z}/3\mathbb{Z}$. 
Let $Y=\{(1,0),(0,1),(0,-1)\}$ 
Let $a=(1,0)$, $b=(1,1)$.  If we choose $H$ to be the copy of $\mathbb{Z}/ 
2\mathbb{Z}$, then the $G/H$ series is $(0,1)$, the $H$ series is $(1,0)$, 
and we can take $((0,1),(1,0))$ as our desired series. 
 
If we take $H$ to be the copy of $\mathbb{Z}/3\mathbb{Z}$, then 
the $G/H$ series is $(1,0)$, and the series for $H$ is $(0,1)$. In 
this case we cannot just concatenate the two series, because we 
are in the undesirable situation described above where $\bar 
b=\bar a$ and the $G/H$ series sums to $a$.  Thus we take the 
first term of the $H$ series (which in this case happens to be all 
of the $H$ series), followed by the $G/H$ series, followed by the 
rest of the $H$ series (which in this case happens to be empty) 
and we obtain $((0,1),(1,0))$ as our desired series. \end{example}

\begin{lemma}\label{spos} The greatest common divisor $s$ of all the $T_i$ 
is greater than 1. \end{lemma} 
 
\begin{proof} 
Suppose otherwise.  Let $i$ be as small as possible, so that 
$t_{i}$ divides $k$. By Lemma ~\ref{t0ne1}, $i>0$.  We will now 
demonstrate that all multiples of $t_{i}$ which are not multiples 
of $t_{i-1}$ must be type 1.  However, since elements of $T_{i}$ 
are of even type, this would force $T_{i}$ to be empty, and 
$t_{i}=t_{i-1}$, a contradiction. 
 
By $\Rb_{i-1}$, adding an element of $T_{i}$ to an element of type 
1 not divisible by $t_{i-1}$ yields another element of type 1. Let 
$x$ be an arbitrary element of $T_{i}$ which is not a multiple of 
$t_{i-1}$. We wish to write $x-k$ as the sum of a series of 
elements from $T_{\leq i-1}$ such that, if the partial sums are 
$z_1,\dots,z_m=x-k$, then for no $l$ is $k+z_l$ divisible by 
$t_{i-1}$.  If we can do this, we can conclude that $x$ is type 1. 
 
We know that the elements of $T_{i}$ generate $t_{i}\mathbb{Z}/t_{i-1}\mathbb{Z}$, 
but in fact more is true.  By Lemma ~\ref{neg}, we know that $T_{i}$ 
contains a set of generators and their negatives for 
$t_{i}\mathbb{Z}/t_{i-1}\mathbb{Z}$.  We can therefore apply Lemma ~\ref{abel}, and 
we are done. \end{proof}

\section{The Cone of $F$-positive Symmetric Functions\label{F-pos}} 
 
We now consider the set $\K$ of all $F \in \Lambda$ having a 
nonnegative representation in terms of the basis of fundamental 
quasisymmetric functions, that is, 
\begin{eqnarray} 
\K = \{ \sum_\alpha c_\alpha F_\alpha \in \Lambda ~|~ c_\alpha \ge 0 
\hbox{\rm ~for all~} \alpha \}. 
\label{Kdef} 
\end{eqnarray} 
 
Since $\K$ is the intersection of $\Lambda$ with the nonnegative 
orthant of $\QQ$ (with respect to the basis $\{F_\beta\}$), $\K_n 
:= \K \cap \Lambda_n$ is a polyhedral cone for each $n\ge 0$. It 
contains the Schur functions $s_\lambda$, $\lambda \vdash n$, so 
it has full dimension in $\Lambda_n$.

\subsection{The generators of $\K_n$.} 
We consider first the minimal generators of the cone $\K_n$, i.e., 
its 1-dimensional faces or extreme rays.  These include all the Schur 
functions and, in general, can be characterized by a condition of 
being balanced. 
 
We begin by considering the notion of the spread of a quasisymmetric function. 
For $\beta \preceq \gamma$, we denote by 
\begin{eqnarray} 
[\beta,\gamma]_\preceq = \{\alpha ~|~ \beta \preceq \alpha \preceq \gamma\} 
\label{lexinterval} 
\end{eqnarray} 
the lexocographic interval between $\beta$ and $\gamma$. For a 
quasisymmetric function $F=\sum c_\alpha F_\alpha \in \QQ$, we 
define the {\it spread} of $F$ to be the smallest lexocographic 
interval $[\beta,\gamma]_\preceq$ so that $c_\alpha =0$ whenever 
$\alpha \notin [\beta,\gamma]_\preceq$. 
 
For a partition $\lambda \vdash n$, we let $\lambda'$ denote 
the conjugate partition and define the composition 
\begin{eqnarray} 
\widetilde{\lambda}:=\beta\left({[n-1]\setminus S\left(\lambda'\right)}\right). 
\label{lambdatilde} 
\end{eqnarray} 
Thus if $\lambda = 33$, then $\lambda'= 222$, so 
${[5]\setminus S(\lambda')}= \{1,3,5\}\subset[5]$ and 
$\widetilde{\lambda}= 1221$.  Note that $\lambda$ corresponds to the 
descent set of the tableaux $T_r$ obtained by filling the Ferrers shape 
$\lambda$ by rows, $\widetilde{\lambda}$ corresponds similarly to 
descent set of the filling $T_c$ by columns and 
$\widetilde{\lambda} \preceq \lambda$, with equality if and only if 
$\lambda$ is $n$ or $1^n$. 
In the example above, we have 
\begin{eqnarray*} 
T_r=\begin{array}{ccc}1&2&3\\4&5&6\\ 
\end{array} 
\hbox{\rm ~and~} 
T_c=\begin{array}{ccc}1&3&5\\2&4&6\\ 
\end{array}. 
\end{eqnarray*} 
Also note that $\widetilde{\lambda} = \widetilde{\nu}$ if and only 
if $\lambda = \nu$.

\begin{proposition} 
The spread of the Schur function $s_\lambda$ is the interval 
$[\widetilde{\lambda},\lambda]_\preceq$. 
\label{spread} 
\end{proposition} 
 
\begin{proof} 
Recall that $s_\lambda = \sum c_\alpha F_\alpha$ where $c_\alpha$ is 
the number of standard Young tableaux $T$ of shape $\lambda$ with 
$\alpha = \beta(D(T))$. 
Let $T_r$, respectively, $T_c$ be the standard Young tableaux obtained by 
filling the Ferrers diagram with shape $\lambda$ by rows, respectively, 
by columns.  As noted, $T_r$ and $T_c$ correspond this way to $\lambda$ and 
$\widetilde{\lambda}$.  Now for any other tableaux $T$, let $i_r$ be the 
first index for which $i_r+1$ is not in the same row as in $T_r$ and 
let $i_c$ be the first index for which $i_c+1$ is not in the same column 
as in $T_c$.  Then $i_r$ is a descent in $T$ but not in $T_r$ and $i_c$ 
is a descent in $T_c$ but not in $T$, so 
$\widetilde{\lambda} \prec \beta(D(T)) \prec \lambda$. 
\end{proof} 
 
\begin{lemma} 
Suppose $\lambda,\nu \vdash n$, $\lambda \ne \nu$, and 
the spread of $s_\nu$ is a subset of the spread of $s_\lambda$. 
Then if $s_\lambda = \sum c_\beta F_\beta$ it follows that $c_\nu = 0$. 
\label{nonminint} 
\end{lemma} 
 
\begin{proof} 
By assumption, we have $\widetilde{\lambda} \prec \widetilde{\nu} 
\preceq \nu \prec \lambda$.  The first inequality implies 
${\nu'} \prec \lambda'$, so there is a minimum 
index $j>1$ so that 
\begin{eqnarray} 
\nu_1 + \cdots + \nu_j > \lambda_1 + \cdots + \lambda_j. 
\label{nuwins} 
\end{eqnarray} 
 
Now if $c_\nu \not= 0$, then there must be a filling $T$ of the shape 
$\lambda$ with $\beta\left(D(T)\right)=\nu$.  Then indices 
$1,2,\dots,\nu_1$ need to be in the first row of $T$, indices 
$\nu_1+1,\dots,\nu_1+\nu_2$ need to be in the first two rows, 
etc.  However (\ref{nuwins}) indicates this filling will fail at 
row $j$. 
\end{proof} 
 
We can now prove the main result of this section. 
 
\begin{theorem} 
The Schur functions $s_\lambda$ are extreme in the cone $\K$. 
\label{schurextreme} 
\end{theorem} 
 
\begin{proof} 
Suppose $s_\lambda = F_1+F_2$ with $F_1,F_2 \in \K$.  Then 
$F_i = \sum_\mu a^i_\mu s_\mu$ with 
\begin{eqnarray} 
a^1_\lambda+a^2_\lambda = 1 \hbox{\rm \quad and\quad} a^1_\mu+a^2_\mu = 0, 
\quad\mu\not=\lambda. 
\label{asums} 
\end{eqnarray} 
Suppose $F_i=\sum c^i_\beta F_\beta$. 
 
If there is a $\mu \not= \lambda$ with $a^i_\mu \not= 0$, then 
either $\mu \succ \lambda$, $\widetilde{\mu} \prec \widetilde{\lambda}$ 
or the spread of $s_\mu$ is a subset of the spread of $s_\lambda$. 
If there is such a $\mu$ with $\mu \succ \lambda$, choose one which is 
lexicographically largest.  If not, but there is one with 
$\widetilde{\mu} \prec \widetilde{\lambda}$, then choose such a $\mu$ such 
that $\widetilde\mu$ is lexicographically smallest.  Otherwise, choose a 
lexicographically largest $\mu$ with 
the spread of $s_\mu$ a subset of the spread of $s_\lambda$. 
By Proposition \ref{spread} and Lemma \ref{nonminint}, 
one of the $F_i$ must have $c^i_\mu < 0$ or $c^i_{\widetilde{\mu}}<0$ 
for the chosen $\mu$. 
 
Thus $a^i_\mu=0$ for $\mu \not= \lambda$ and so 
both $F_1$ and $F_2$ are multiples of $s_\lambda$, showing $s_\lambda$ 
to be extreme. 
\end{proof} 
 
Note that there are extremes other than the Schur functions.  The first one appears
when $n=4$:
$$s_{31} + s_{211} - s_{22} = F_{31} + F_{13} + F_{211} + F_{112}$$
is extreme in $\K_{4}$.  In $\K_{5}$, there are two such extremes,
$s_{311} + s_{2111} - s_{221}$ and $s_{41} + s_{311} - s_{32}$.  In $\K_{6}$,
there are 23.  At present there is no general description of which combinations of Schur
functions are extreme. 
 
Consequently, we consider next the problem of determining when a quasisymmetric 
function $F=\sum h_S F_S$ is an extreme element of the cone $\K$ 
of $F$-positive symmetric functions.  (Here we begin indexing by 
subsets of $[n]$  in place of compositions of $n+1$, 
where $F_S=F_{\beta(S)}$.)  We relate this to a property 
of the multicollection $\{S^{~h_S} \}$, which leads to the notion 
of {\em fully balanced multicollections} of subsets of a finite 
set. Fully balanced multicollections with nonnegative 
multiplicities will yield $F$-positive symmetric functions, in 
general, while minimal such collections give rise to extremes. 
 
We say a subset $S\subset [n]$ has {\em profile} $a_1,\dots,a_k$ if $S$ consists 
of maximal consecutive strings of length $a_1,\dots, a_k$ in some order. 
In this case, $|S| = a_1 + \cdots +a_k$.  For example, 
$\{2,3,5,7,8,9\}\subset [11]$ has profile $321000$. 
For $\lambda=\lambda_1\lambda_2\dots\lambda_k \vdash n+1$, define 
\begin{eqnarray} 
\F_\lambda = \{ S\subset [n] ~|~ S \hbox{\rm~ has profile ~} \lambda_1-1,\dots,\lambda_k-1\}. 
\label{F-lambda} 
\end{eqnarray} 
Thus if $S=\{2,3,5,7,8,9\}\subset [11]$, then $S \in \F_{432111}$. 
Further $\F_{11\dots 1}=\{\emptyset \}$,  $S\in \F_{21\dots 1}$ if and only if $|S|=1$, 
and $S\in \F_{221\dots 1}$ if and only if $S=\{i,j\}$, where $i<j-1$, 
while $S=\{ i,i+1 \} \in \F_{31\dots 1}$. 
 
We denote a {\em multicollection} of subsets of $[n]$ by 
$\{S^{~k_S}~|~ S\subset [n] \}=\{ S^{~k_S} \}$, where $k_S$ 
denotes the multiplicity of the subset $S$. 
For our purposes, a 
multicollection $\{S^{~k_S}\}$ can have any rational multiplicities $k_S$. 
\begin{definition} Let $\lambda \vdash n+1$. 
 A multicollection $\{S^{~k_S}\}$  of subsets of $[n]$ is 
\emph{$\lambda$-balanced} if there is a constant $\kappa_\lambda$ such that for all 
$T\in \F_\lambda$, 
\begin{eqnarray} 
\sum_{S\supseteq T} k_S = \kappa_\lambda. 
\label{lambda-sum} 
\end{eqnarray} 
The multicollection $\{S^{~k_S}\}$ is \emph{fully balanced} if it 
$\lambda$-balanced for all $\lambda \vdash n+1$. 
\label{lambda-bal} 
\end{definition}

Multicollections that are $21\dots 1$-balanced have been called \emph{balanced} 
in the literature of cooperative game theory \cite{YK}, although there the term is applied 
to the underlying collection whenever positive multiplicities $k_S$ exist. 
 
\begin{theorem} 
A homogeneous quasisymmetric function $F=\sum_S h_S F_S \in \QQ_{n+1}$ 
is symmetric if and only 
if the multicollection $\{ S^{~h_S}\}$ of subsets of $[n]$ is fully balanced. 
\end{theorem} 
 
\begin{proof} 
Note that, for $\mu \vdash n+1$, $R\in \F_\mu$ if and only if 
$\lambda \left(\beta([n]\setminus R) \right) = \mu$. 
Further, note that if 
$T\in \F_\lambda$ and $R\subset T$, $R\not=T$, then $R\in \F_\mu$ for some 
$\mu \prec \lambda$. 
By inclusion-exclusion, we get 
\begin{eqnarray} 
\sum_{S \supseteq T} h_S = \sum_{R \subseteq T} (-1)^{|R|} 
\sum_{S \subseteq [n] \setminus R} h_S 
= \sum_{R \subseteq T} (-1)^{|R|} f_{[n] \setminus R}, 
\label{hequation} 
\end{eqnarray} 
where $f_S$ and $h_S$ are related 
as $d_\beta$ and $c_\beta$ in (\ref{inversion}).  Now, $F $ is symmetric if and only if 
$f_{[n]\setminus R}$ only depends on $\mu$ for $R\in \F_\mu$. 
Thus if $F$ is symmetric, then (\ref{hequation}) shows the 
sum $\sum_{S \supseteq T} h_S$ to depend only on $\lambda$ (and $\mu \prec \lambda$) 
when $T\in \F_\lambda$. 
 
Now suppose the multicollection $\{S^{~h_S}\}$ is $\lambda$-balanced for all 
$\lambda \vdash n+1$.  We argue by induction on the lexicographic order on 
partitions.  We assume $f_{[n]\setminus R}$ only depends on $\mu$ for all 
$R\in \F_\mu$, $\mu \prec \lambda$.  (The base case for $\lambda=11\dots1$ is 
trivial.)  For $T \in \F_\lambda$, the assertion now follows from (\ref{hequation}), since 
the number of $R \subset T$ with $R\in \F_\mu$, for $\mu \prec \lambda$, depends only on 
$\lambda$. 
\end{proof} 
 
Thus, elements of $\K_{n+1}$ correspond to fully balanced collections with 
nonnegative multiplicities.  Those  with minimal support $\{S~|~ h_S\not= 0\}$ 
correspond to the extremes of the cone.  One can view integral extremes of 
$\K_{n+1}$ as combinatorial designs of an extremely balanced sort: each element 
of $[n]$ is in the same number of sets (counting multiplicity), as are each nonadjacent 
pair, each adjacent pair, etc.  One is led to wonder whether the designs coming this 
way from Schur functions have special properties among these.  The first of these 
for which the multiplicities are not all one is 
\begin{eqnarray*} 
s_{321}&=&F_{\{1,3\}}+F_{\{1,4\}}    + F_{\{2,3\}}+ 2F_{\{2,4\}}+ F_{\{2,5\}} +F_{\{3,4\}} 
 +F_{\{3,5\}}  \\ 
  &+&F_{\{1,2,4\}} +F_{\{1,2,5\}} +F_{\{1,3,4\}} 
 +2F_{\{1,3,5\}} +F_{\{1,4,5\}} +F_{\{2,3,5\}}  +F_{\{2,4,5\}}. 
 \end{eqnarray*} 
 Here $\kappa_{21111}=8$, $\kappa_{3111}=2$, $\kappa_{2211}=4$, $\kappa_{321}=1$ 
 and $\kappa_{222}=2$. 
 
\subsection{The facets of $\K_n$.} 
To describe the facets of $\K_n$, we rewrite (\ref{Kdef}) as follows. 
Since the Schur functions $s_\lambda$, $\lambda \vdash n$, are a basis 
for $\Lambda_n$, writing 
$s_\lambda = \sum_\beta [s_\lambda]_{F_\beta} F_\beta$, 
we see that 
\begin{eqnarray} 
\K_n = \{ \sum_{\lambda \vdash n} c_\lambda s_\lambda ~|~ 
        \sum_\lambda c_\lambda [s_\lambda]_{F_\beta} \ge 0 
            \hbox{\rm ~for all~} \beta \vDash n\}. 
\label{Kineq} 
\end{eqnarray} 
Equation (\ref{Kineq}) gives $2^{n-1}$ inequalities for $\K_n$, one 
for each $\beta \vDash n$.  However, when $\beta \sim \gamma$, these 
inequalities are identical (see Definition \ref{beta-equiv}).  In fact, 
we conjecture that these are the only redundant inequalities, so the 
facets of $\K_n$ would be in bijection with the equivalence classes of 
compositions under $\sim$.

The inequality for $\K_n$ given by $c_\alpha \ge 0$ in (\ref{Kdef}) 
is redundant if and only if there exist $a_\beta \ge 0$ such that 
\begin{eqnarray} 
c_\alpha = \sum_{\beta \not\sim \alpha} a_\beta c_\beta 
\label{redundant} 
\end{eqnarray} 
holds for all 
$F=\sum c_\gamma F_\gamma \in \Lambda$. 
 
For each composition $\beta \vDash n$ we define the vector 
$v_\beta = (v_{\beta,\lambda}; \lambda \vdash n)$ by 
$v_{\beta,\lambda} = \mult_{\M(\beta)}(\lambda)$.  By definition, 
$\beta \sim \gamma$ if and only if $v_\beta = v_\gamma$. 
 
\begin{proposition} 
The inequality $c_\alpha\ge 0$ is redundant for some $\alpha 
\vDash n$ if and only if $v_\alpha$ is not extreme in the convex 
hull of all $v_\beta$, $\beta \vDash n$. \label{vertex} 
\end{proposition} 
 
We end with the following conjecture. 
 
\begin{conjecture} 
\label{facets-conj} 
Any one, and so all, of the equivalent statements holds: 
\begin{enumerate} 
\item 
The facets of $\K_n$ are in bijection with the equivalence classes 
of compositions $\beta \vDash n$, 
\item 
The inequalities $c_\alpha\ge 0$, $\alpha \vDash n$, are all irredundant, 
\item 
Each $v_\alpha$ is extreme in the convex hull of all $v_\beta$, $\beta 
\vDash n$. 
\end{enumerate} 
\end{conjecture} 
 
One can imagine an approach to 
Conjecture \ref{facets-conj} that uses Theorem  \ref{decomp-thm} 
along with a separation argument for $v_\beta$ 
that targets the decomposition structure of the composition $\beta$. 
 
\subsection{Ribbon-positivity} 
 
One can call a symmetric function $F$-\emph{positive} if it belongs to the cone $\K$. 
Being $F$-positive is a weakening of the condition of being Schur-positive. 
Both $F$-positivity and Schur-positivity are closed under taking products.
The quasisymmetric functions $F_S$ have been identified with characters
of 0-Hecke algebras \cite[\S 4.1]{Thibon}, lending a
representation-theoretic
interpretation to being $F$-positive.

A strengthening of Schur-positivity would be 
what we might call \emph{ribbon-positivity}, that is, belonging to the cone 
in $\Lambda$ spanned by the ribbon Schur functions $r_{\alpha}, \alpha \vDash n$. 
By \eqref{ribbonmult}, this cone is also closed under taking products. 
As far as we know, this ribbon cone has not been studied. 
 
Stronger still would be 
membership in the cone spanned by the $r_{\lambda}, \lambda \vdash n$. 
That this simplicial cone is strictly smaller than the ribbon cone can be seen 
from the fact that $r_{132} = r_{321} +r_{33} - r_{42}$, which follows from \eqref{ribbonlinear}. 
This relation suggests that the ribbon Schur functions might all be extreme in the ribbon cone. 
 
We know of one example of ribbon-positive symmetric functions.  In \cite{Hersh}, Hersh shows 
that the chain-enumeration quasisymmetric function of $k$-shuffle posets are ribbon-positive 
by showing they are sums of products of ribbon Schur functions. 
 
\providecommand{\bysame}{\leavevmode\hbox 
to3em{\hrulefill}\thinspace}


\begin{thebibliography}{10} 
 
\bibitem{BLiu} L.J.\ Billera and N. Liu, Noncommutative enumeration in 
graded posets, {\it J. Algebraic Combin.} {\bf 12} (2000), 7--24. 
 
\bibitem{Ful} W.\ Fulton, {\it Young Tableaux},  Cambridge University 
Press, Cambridge, UK, 1997. 
 
\bibitem{Gap} 
The GAP~Group, \emph{GAP -- Groups, Algorithms, and Programming, 
Version 4.3}; 2002 \verb+(http://www.gap-system.org)+. 
 
\bibitem{GKL} I.\ M.\ Gel'fand, D.\ Krob, A.\ Lascoux, B.\ Leclerc, 
V.\ Retakh and J.-Y.\ Thibon, Noncommutative symmetric functions, 
{\it Adv. Math.} {\bf 112} (1995), 218--348. 
 
 
\bibitem{GR} I.\ Gessel and C.\ Reutenauer, Counting permutations 
with given cycle structure and descent set, {\it J. Combin. 
Theory Ser.\ A} {\bf 64} (1993), 189--215. 
 
\bibitem{Hersh} P.\ Hersh, Chain decomposition and the flag $f$-vector, 
{\it J. Combin. Theory Ser.\ A} {\bf 103} (2003), 27--52.  
 
\bibitem{YK} Y.\ Kannai, The core and balancedness, Chap. 12 in {\it Handbook of Game 
Theory with Economic Applications}, Vol. 1, R.J.\ Aumann and 
S.\ Hart, eds., Elsevier Science Publishers (North-Holland), Amsterdam, 1992. 
 
\bibitem{KKN} A.N.\ Kirillov. A.\ Kuniba and T.\ Nakanishi, Skew diagram method 
in spectral decomposition of integrable lattice models, {\it Comm.\ Math.\ Phys.\ } 
{\bf 185} (1997), 441-465. 
 
%
 
\bibitem{LP} A.\ Lascoux and P.\  Pragacz,  Ribbon Schur functions, 
{\it European J. Combin.} {\bf 9} (1988), no. 6, 561--574. 
 
\bibitem{MacD} I.\ Macdonald, {\it Symmetric Functions and Hall Polynomials, 2nd Edition}, 
Oxford University Press, New York, USA, 1995. 
 
\bibitem{MM} P.A.\ MacMahon, {\it Combinatory Analysis}, Cambridge University Press, 1917 
(Vol. I), 1918 (Vol. II), reprint, Chelsea, New York, USA, 1960. 
 
 
\bibitem{ECII} R.\ Stanley, {\it Enumerative Combinatorics, Vol. 2}, 
Cambridge Studies in Advanced Mathematics, Vol. 62, Cambridge 
University Press, Cambridge, UK, 1999. 
 
 
 
\bibitem{Thibon} J.-Y.\ Thibon, Lectures on noncommutative symmetric
functions,
{\it Interaction of combinatorics and representation theory},  39--94,
MSJ Mem. 11,
Math. Soc. Japan, Tokyo, 2001.


\bibitem{SvW} S.\ van Willigenburg, Equality of Schur and skew Schur functions, 
{\it Ann. Comb.} (to appear). 


 
 
 
\end{thebibliography}
\end{document}